\documentclass[11 pt]{amsart}
\usepackage{amsmath}
\usepackage{amsfonts}
\usepackage{enumerate}
\usepackage{verbatim}
\usepackage[pdftex]{graphicx}
\DeclareGraphicsExtensions{.png,}
\newtheorem{lemma}{Lemma}
\newtheorem{theorem}{Theorem}
\newtheorem{example}{Example}
\newtheorem{prop}{Proposition}
\newtheorem{corollary}{Corollary}

\newtheorem{defn}{Definition}
\numberwithin{equation}{section}
\numberwithin{defn}{section}
\numberwithin{prop}{section}
\numberwithin{lemma}{section}

\newcommand{\R}{\mathbb{R}}

\newcommand{\N}{\mathbb{N}}

\newcommand{\M}{\mathcal{M}}
\newcommand{\T}{\mathbb{T}}
\newcommand{\var}{\text{var}}
\newcommand{\D}{\mathcal{D}}
\newcommand{\W}{\mathcal{W}}

\newcommand{\A}{\mathcal{A}}
\newcommand{\V}{\mathcal{V}}

\newcommand{\B}{\mathcal{B}}

\newcommand{\E}{\mathcal{E}}

\renewcommand{\dim}{\mbox{dim}_{\mathcal{H}}}
\newcommand{\Dim}{\mbox{dim}_{\mathcal{P}}}
\newcommand{\BDim}{\overline{\mbox{dim}}_{B}}

\title[The Packing Spectrum for Birkhoff Averages]{The Packing Spectrum for Birkhoff Averages on a self-affine Repeller}
\author{Henry WJ Reeve}
\address{Henry WJ Reeve\\Department of Mathematics\\ The University of Bristol\\
University Walk\\Clifton\\ Bristol\\BS8 1TW\\UK.}
\email{henrywjreeve@googlemail.com}

\thanks{The author would like to thank Thomas Jordan and Micha\l{} Rams for all their help in preparing this paper. This project was started in the Instytut Matematyczny PAN and I would also like to thank the Institute, especially Feliks Przytycki, for their kind hospitality. I would also like to thank the Engineering and Physical Sciences Research Council and the Conformal Structures and Dynamics network for their financial support. In addition I would like to thank the referee, whose careful reading of this paper has greatly improved the presentation.}

\begin{document}
\begin{abstract} We consider the multifractal analysis for Birkhoff averages of continuous potentials on a self-affine Sierpi\'{n}ski sponge. In particular, we give a variational principal for the packing dimension of the level sets. Furthermore, we prove that the packing spectrum is concave and continuous. We give a sufficient condition for the packing spectrum to be real analytic, but show that for general H\"{o}lder continuous potentials, this need not be the case. We also give a precise criterion for when the packing spectrum attains the full packing dimension of the repeller. Again, we present an example showing that this is not always the case. 
\end{abstract}
\maketitle
\section{Introduction and Statement of Results}
Let $\Lambda$ be the repeller of a $C^{1+\epsilon}$ map $f: X \rightarrow X$. Given some continuous potential $\varphi: \Lambda \rightarrow \R^N$ and some $\alpha \in \R^N$ we are interested in the set of points in the repeller for which the Birkhoff average converges to $\alpha$,
\begin{equation}
E_{\varphi}(\alpha):= \bigg \{ x \in \Lambda : \lim_{n \rightarrow \infty} \frac{1}{n} \sum_{q=0}^{n-1} \varphi(f^q(x))=\alpha\bigg\}
.\end{equation}
We would like to understand how the geometric complexity of $E_{\varphi}(\alpha)$ varies as a function of $\alpha$. Geometric complexity, here, is to be understood in terms of the dual notions of Hausdorff dimension $\dim$, defined in terms of minimal coverings, and packing dimension $\Dim$, defined in terms of maximal packings (see \cite[ Section 6.2]{Edgar1} or \cite[Chapter 3]{Falconer}). We refer to $\alpha\mapsto \dim E_{\varphi}(\alpha)$ as the Hausdorff spectrum and $\alpha \mapsto \Dim E_{\varphi}(\alpha)$ as the packing spectrum.

In the conformal setting there is a well known variational principle giving the values for both spectra. To recall this result we require some terminology.
Let $\Lambda$ be a repeller for an expanding $C^{1+\epsilon}$ map $f$ of a smooth manifold $X$. We let $\M(\Lambda,f)$ denote the set of $f$-invariant Borel probability measures supported on $\Lambda$. Given $\mu \in \M(\Lambda,f)$ we let $h_{\mu}(f)$ denote the Kolmogorov-Sinai entropy (see \cite[Section 4.10]{Walters}) and let $\lambda_{\mu}(f):=\int \log ||f'|| d\mu$ denote the Lyapunov exponent. Given a continuous potential $\varphi:\Lambda\rightarrow \R$ we let $A(\varphi):=\{\int \varphi d\mu:\mu\in\M(\Lambda,f)\}$. We may now recall the classic result due to Pesin and Weiss \cite{Pesin Weiss Birkhoff}, Fan, Feng and Wu \cite{Fan Feng Wu}, Barreira and Saussol \cite{Barreira Saussol}, Feng, Lau and Wu \cite{Feng Lau Wu} and Olsen \cite{Olsen Multifractal 1}, \cite{Olsen Multifractal 4}. See also \cite{non-uniformly hyperbolic} and \cite{Gelfert Rams} for extensions to classes of non-uniformly hyperbolic functions.
\begin{theorem}[Feng, Lau, Wu]\label{Conformal Birkhoff}
Let $\Lambda$ be the repeller for an expanding $C^{1+\epsilon}$ map $f:X\rightarrow X$. Suppose that $f$ is conformal and topologically mixing on $\Lambda$. Then for all $\alpha\in A(\varphi)$ we have
\begin{eqnarray*}
\dim E_{\varphi}(\alpha)= \Dim E_{\varphi}(\alpha)= \sup \left\lbrace \frac{h_{\mu}(f)}{\lambda_{\mu}(f)}: \mu \in \M(\Lambda,f), \int \varphi d\mu =\alpha \right\rbrace.
\end{eqnarray*}
\end{theorem}
In particular, when $f$ is conformal and uniformly hyperbolic the Hausdorff and packing spectra coincide. This is a consequence of the fact that the level set $E_{\varphi}(\alpha)$ corresponds to a type of statistical convergence, together with the neat relationship between geometric and statistical properties which holds in the conformal setting. By contrast, the packing and Hausdorff dimensions of level sets defined in terms of divergent asymptotic properties may differ (see \cite{Snigreva Packing} and \cite{Olsen Multifractal 4}). Theorem \ref{Conformal Birkhoff} also allows one to deduce various regularity properties of the spectrum (\cite{Pesin Weiss MA Equilibrium}, \cite{Fan Feng Wu}, \cite{Barreira Saussol}, \cite{Olsen Multifractal 1}). The spectrum is continuous and when $\varphi$ is H\"{o}lder continuous it is also real analytic. When the Lyapunov exponent $\lambda_{\mu}(f)$ is given by a fixed constant, independent of $\mu \in \M(\Lambda,f)$, the spectrum is also concave.

The dimension theory of non-conformal systems, for which there is no simple correspondence between geometric and statistical properties, is much less well understood. For the most part we only have almost all type results, both for the dimension of a repeller  \cite{Falconer singular value} and for the dimension of level sets for Birkhoff averages \cite{Jordan Simon}. One class of non conformal fractals for which we do have deterministic results (\cite{Bedford}, \cite{King}, \cite{McMullen}, \cite{Jordan Rams}, \cite{Barral Mensi}, \cite{Nielsen}, \cite{Kenyon Peres}) are the self-affine Sierpi\'{n}ski sponges introduced by Bedford \cite{Bedford} and McMullen \cite{McMullen} and generalized by Kenyon and Peres \cite{Kenyon Peres}.
\begin{defn}[Self-affine Sierpi\'{n}ski sponges]\label{Sierpinski sponge}
Let $\T^d:=\R^d/\mathbb{Z}^d$ denote the $d$ dimensional torus. Choose natural numbers $a_1 > a_2 > \cdots > a_d \geq 2$. Let $f$ denote the integer valued diagonal map given by 
\begin{equation} (x_q)_{q=1}^d\mapsto (a_qx_q)_{q=1}^d \text{ for } (x_q)_{q=1}^d \in \T^d
.\end{equation}
Given a digit set $\D \subseteq \prod_{q=1}^d \left\lbrace0,\cdots,a_q-1\right\rbrace$
there is a corresponding self-affine repeller $\Lambda$ given by 
\begin{equation}
\Lambda:=\left\lbrace \left(\sum_{\nu=1}^{\infty}\frac{i_q(\nu)}{a_q^{\nu}}\right)_{q=1}^d:(i_q(\nu))_{q=1}^d\in\D \text{ for all }\nu\in\N \right\rbrace
.\end{equation}
A limit set $\Lambda$ defined in this way is referred to as a self-affine Sierpi\'{n}ski sponge. A two dimensional Sierpi\'{n}ski sponge is known as a Bedford-McMullen carpet.
\end{defn}
\begin{figure}
\includegraphics[width=140mm]{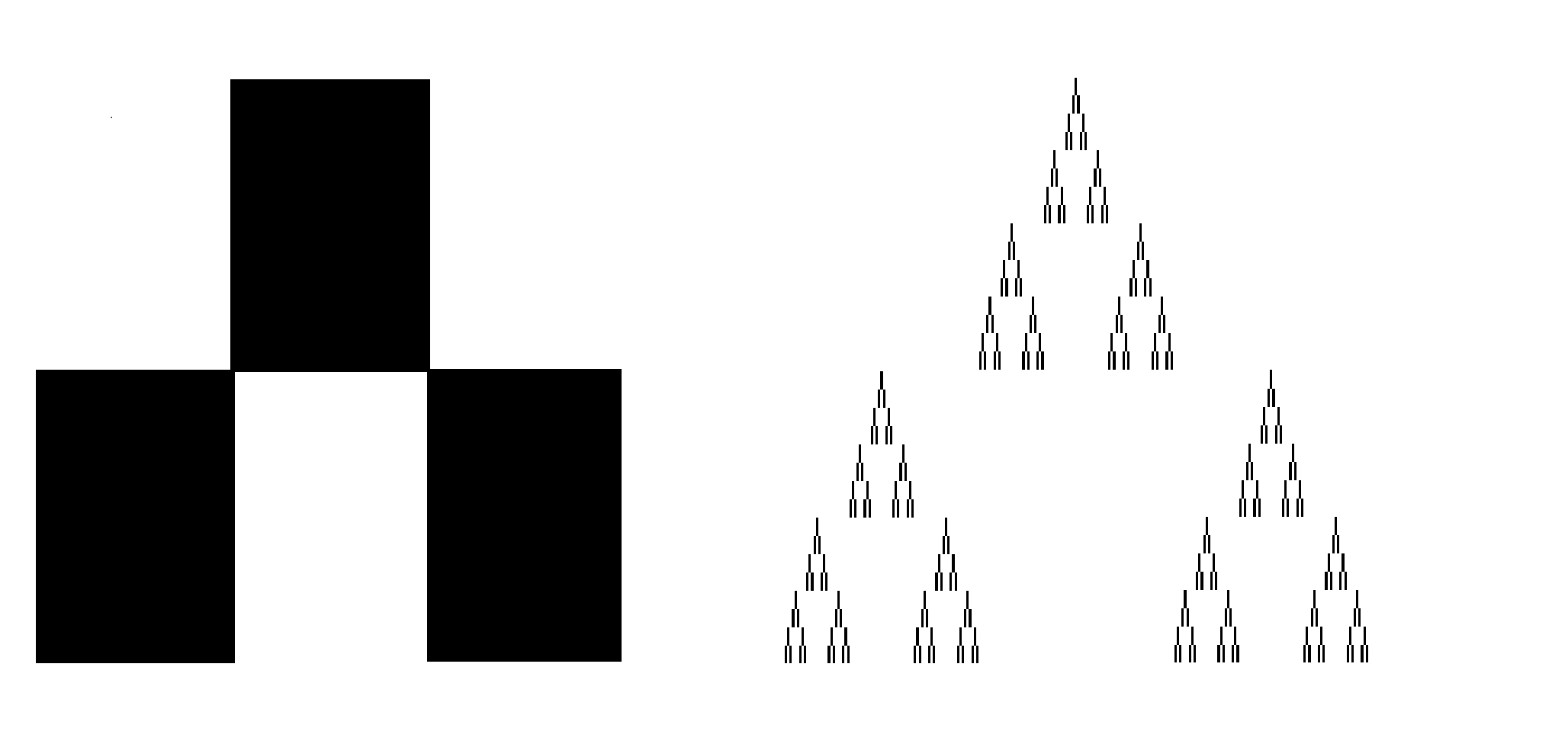}
\caption{A representation of a digit set (left) and the corresponding Bedford-McMullen carpet (right).}
\end{figure}
To state the relevant results concerning self-affine Sierpi\'{n}ski sponges we first introduce some terminology. 
Given a continuous transformation $T$ of a metric space $X$ we let $\M(X,T)$ denote the set of $T$-invariant Borel probability measures, and let $\E(X,T)$ denote the set of $\mu\in \M(X,T)$ which are ergodic.
 
Given $k \leq d$ we let $\pi_k:\T^{d}\rightarrow \T^{d-(k-1)}$ be the projection $\pi_k:(x_q)_{q=1}^d\mapsto (x_q)_{q=k}^d$. We let $f_k:\T^{d-(q-1)}\rightarrow \T^{d-(q-1)}$ denote the map $(x_q)_{q=k}^d \mapsto(a_q x_q)_{q= k}^d$. 

Suppose $A \subseteq \M(\Lambda,f)$. For each $k \leq d$ we define $H^k(f,A)$ by  
\begin{eqnarray*}
H^k(f,A)&:=&\sup \left\lbrace h_{\mu\circ \pi_k^{-1}}(f_k): \mu \in A \right\rbrace.
\end{eqnarray*}
When $A= \M(\Lambda,f)$ we write $H^k(f,A)=H^k(f)$.
 
The following result is due to Kenyon and Peres \cite{Kenyon Peres}.
\begin{theorem}[Kenyon, Peres] \label{KP result} Let $\Lambda$ be a self-affine Sierpi\'{n}ski sponge. Then,
\begin{eqnarray*}
\dim(\Lambda)&=&\sup_{\mu \in \M(\T^d,f)}\left\lbrace \frac{h_{\mu}(f)}{\log a_1}+ \sum_{k=2}^d \left(\frac{1}{\log a_{k}}-\frac{1}{\log a_{k-1}}\right)h_{\mu \circ \pi_k^{-1}}(f_k)\right\rbrace,
\\
\Dim(\Lambda)&=&\frac{H^1(f)}{\log a_1}+ \sum_{k=2}^d \left(\frac{1}{\log a_{k}}-\frac{1}{\log a_{k-1}}\right)H^k(f)
.\end{eqnarray*}
\end{theorem}
Bedford \cite{Bedford} and McMullen \cite{McMullen} independently determined both the Hausdorff dimension and the upper box dimension in the two dimensional setting. In \cite{Kenyon Peres} Kenyon and Peres extend these results to higher dimensions. It follows from \cite[Proposition 3.6]{Falconer Techniques} that the formula for upper box dimension also gives an expression for the packing dimension.

The multifractal analysis of Birkhoff averages is closely related to the multifractal analysis of pointiwise dimension. Given an invariant measure $\nu\in \M(\Lambda,f)$ on a self-affine Sierpi\'{n}ski sponge $\Lambda$ the Hausdorff and packing spectrums for pointwise dimension are given by $\alpha \mapsto \dim D_{\nu}(\alpha)$ and $\alpha \mapsto \Dim D_{\nu}(\alpha)$, respectively, where
\begin{equation}
D_{\nu}(\alpha):= \left\lbrace x \in \Lambda: \lim_{r\rightarrow 0} \frac{\log \nu(B(x,r))}{\log r}=\alpha \right\rbrace.
\end{equation}
In \cite{King} King determined the Hausdorff spectrum for Bernoulli measures on a Bedford-McMullen carpet with strong separation conditions. Olsen extended King's result to Bernoulli measures on $d$ dimensional self-affine Sierpi\'{n}ski sponge \cite{Olsen Affine Pointwise}. The Hausdorff spectrum for Gibbs measures was determined by Barral and Mensi \cite{Barral Mensi Pointwise} for Bedford-McMullen carpets, and by Barral and Feng \cite{Barral Feng} for a $d$ dimensional self-affine Sierpi\'{n}ski sponge. In \cite{Jordan Rams} Jordan and Rams gave the Hausdorff spectrum for Bernoulli measures on a Bedford-McMullen carpet without the strong separation conditions required in \cite{King}, \cite{Olsen Affine Pointwise}, \cite{Barral Mensi Pointwise} and \cite{Barral Feng}. In contrast almost nothing is known about the packing spectrum for pointwise dimension on a self-affine Sierpi\'{n}ski sponge. In this paper we determine $\alpha \mapsto \Dim D_{\nu}(\alpha)$ for a very limited class of Bernoulli measures $\nu$ on self-affine Sierpi\'{n}ski sponges with strong separation conditions. We also give an example disproving a conjecture of Olsen \cite[Conjecture 4.1.7]{Olsen Affine Pointwise} (see Section \ref{Generalisations}). However, the main focus for this article is the packing spectrum for Birkhoff averages.
 
The first result concerning the multifractal analysis of Birkhoff averages for self-affine Sierpi\'{n}ski sponges is due to Nielsen \cite{Nielsen}. Suppose $\Lambda$ is a self-affine Sierpi\'{n}ski sponge. For $x\in\Lambda$ we let 
\begin{equation*}
\Gamma(x):=\left\lbrace ((i_q(\nu))_{q=1}^d)_{\nu\in\N}: x= \left(\sum_{\nu=1}^{\infty}\frac{i_q(\nu)}{a_q^{\nu}} \right)_{q=1}^d \right\rbrace
.\end{equation*}
Given a probability vector ${\mathbf{p}}=(p_{l})_{l\in\D}$ defined over a digit set $\D$ we let $N_{l}\left(\omega|n\right):=\#\{\nu\leq n:\omega_{\nu}=l\}$, where $\#$ denotes cardinality, and define 
\begin{eqnarray*}
\Lambda_{\mathbf{p}}:=\left\lbrace x\in\Lambda:\exists \omega\in \Gamma(x)\text{ with } \lim_{n\rightarrow \infty} \frac{N_{l}\left(\omega|n\right)}{n}=p_{l} \text{ for each } l\in\D \right\rbrace
.\end{eqnarray*} Let $\mu_{\mathbf{p}}$ denote the Bernoulli measure on $\Lambda$ corresponding to the probability vector ${\mathbf{p}}$.
In \cite{Nielsen} Nielsen proved the following formula for the Hausdorff and packing dimension of $\Lambda_{\mathbf{p}}$ in the two dimensional case. With minor modifications the proof also applies in higher dimensions.
\begin{theorem}[Nielsen] \label{Nielsen result}
\begin{eqnarray*}
\Dim \Lambda_{\mathbf{p}}=\dim \Lambda_{\mathbf{p}}
= \frac{h_{\mu_{\mathbf{p}}}(f)}{\log a_1} +\sum_{k=2}^d\bigg(\frac{1}{\log a_k}-\frac{1}{\log a_{k-1}}\bigg) h_{\mu_{\mathbf{p}}\circ \pi_k}(f_k).
\end{eqnarray*}
\end{theorem}
In particular, this shows that for a certain special class of Birkhoff averages, defined over a self-affine Sierpi\'{n}ski sponge, we always have $\Dim E_{\varphi}(\alpha) = \dim E_{\varphi}(\alpha)$. However, it follows from Theorem \ref{KP result} that for self-affine Sierpi\'{n}ski sponges we often have $\dim(\Lambda)<\Dim(\Lambda)$. Consequently, by considering any $\varphi$ which is cohomologous to a constant, the Hausdorff and packing spectra for Birkhoff averages on a Bedford-McMullen repeller do not always coincide. This is a consequence of there being two distinct rates of expansion. It takes less time for a difference along the direction of strong repulsion to be blown up to the scale of the Markov partition than it does for a similarly sized difference along the direction of weak repulsion. As such a given geometric scale will correspond to two time scales, often resulting in a difference between Hausdorff and packing dimensions. The reason that this does not affect the coincidence of $\Dim \Lambda_{\mathbf{p}}$ and $\dim \Lambda_{\mathbf{p}}$ is that the convergence of $N_d(\omega|n)/n$ forces points in the set $\Lambda_{\mathbf{p}}$ to display similar behaviour at both time scales. For less restrictive level sets this need not be the case. 

This extra level of complexity in the non conformal case makes the question of Hausdorff and packing spectra for Birkhoff averages on Bedford-McMullen repellers an interesting one, where we do not expect to observe the same behaviour as in the conformal case. The first part of this question was answered by Barral, Feng and Mensi in \cite{Barral Mensi} and \cite{Barral Feng}. Given an integer valued diagonal map $f$ on a self-affine Sierpi\'{n}ski sponge $\Lambda$ and a continuous potential $\varphi:\Lambda\rightarrow \R^N$ we let $A(\varphi):=\{\int \varphi d\mu:\mu\in\M(\Lambda,f)\}$. One can easily see that $E_{\varphi}(\alpha)=\emptyset$ for $\alpha\notin A(\varphi)$. The following result concerning the Hausdorff spectrum is due to Barral and Feng \cite{Barral Feng}. 
\begin{theorem}[Barral, Feng]\label{Barral Feng}
 Let $\Lambda$ be a self-affine Sierpi\'{n}ski sponge. Let $\varphi:\Lambda \rightarrow \R^N$ be a continuous potential. Then for all $\alpha \in A(\varphi)$ we have
\begin{equation*}
\dim E_{\varphi}(\alpha)=\sup\left\lbrace \frac{h_{\mu}(f)}{\log a_1}+ \sum_{k=2}^d \left(\frac{1}{\log a_{k}}-\frac{1}{\log a_{k-1}}\right)h_{\mu \circ \pi_k^{-1}}(f_k) \right\rbrace,\\
\end{equation*}
where the supremum is taken over all $ \mu\in\M(\Lambda,f)$ with $\int \varphi d\mu=\alpha$.
\end{theorem}
This extends the work of Barral and Mensi in \cite{Barral Mensi} where the Hausdorff spectrum for H\"{o}lder continuous potentials on a Bedford-McMullen carpet is given as the Legendre transform of an explicit moment function.  

In this paper we prove a dual result for the packing spectrum.
For each $\alpha \in A(\varphi)$ we define $H^k(T,\varphi,\alpha)$ for $k=1,\cdots, d$ by  
\begin{eqnarray*}
H^k(f,\varphi,\alpha)&:=&\sup \left\lbrace h_{\mu\circ \pi_k^{-1}}(f_k): \mu \in \M_{T}(\Lambda), \int \varphi d\mu =\alpha \right\rbrace.
\end{eqnarray*}
\begin{theorem}\label{Packing level sets}
Let $\Lambda$ be a self-affine Sierpi\'{n}ski sponge. Let $\varphi:\Lambda \rightarrow \R^N$ be some continuous potential. Then for all $\alpha \in A(\varphi)$ we have
\[\Dim{E_{\varphi}(\alpha)}=  \frac{H^1(f,\varphi,\alpha)}{\log a_1} +\sum_{k=2}^d \left(\frac{1}{\log a_k}-\frac{1}{\log a_{k-1}}\right) H^k(f,\varphi,\alpha).\]
\end{theorem}
In fact Theorem \ref{Packing level sets} follows from the more general Theorem \ref{Packing Main}. Given a Borel probability measure $\mu \in \M(\Lambda)$ we define 
\begin{equation}
A_n(\mu):= \frac{1}{n}\sum_{k=0}^{n-1} \mu \circ f^{-k}.
\end{equation}
Given $x \in \Lambda$ we let $\V(x)$ denote the set of all weak $*$ accumulation points of the sequence of measures $(A_n(\delta_x))_{n\in\N}$ where $\delta_x$ denotes the Dirac measure concentrated at $x$. Note that $\V(x) \subseteq \M(\Lambda,f)$ \cite[Theorem 6.9]{Walters} for all $x \in \Lambda$. Given $A \subseteq \M(\Lambda,f)$ we define 
\begin{eqnarray}
X(A):&=& \left\lbrace x \in \Lambda : \V(x)=A\right\rbrace, \\
Y(A):&=& \left\lbrace x \in \Lambda : \V(x)\subseteq A\right\rbrace \nonumber.
\end{eqnarray}
In \cite{Barral Feng} Barral and Feng considered the special case in which $A=\{\mu\}$ for some $\mu \in \M(\Lambda,f)$. It follows that $X(\{\mu\})=Y(\{\mu\})$ and the Hausdorff and packing dimensions coincide. However, in general this is not the case.

\begin{theorem}\label{Packing Main}
Let $\Lambda$ be a self-affine Sierpi\'{n}ski sponge. Suppose that $A$ is a non-empty closed convex subset of $\M(\Lambda,f)$. Then, 
\[\Dim{X(A)}=\Dim Y(A)= \frac{H^1(f,A)}{\log a_1} +\sum_{k=2}^d \left(\frac{1}{\log a_k}-\frac{1}{\log a_{k-1}}\right) H^k(f,A).\]
\end{theorem}

The central difficulty in determining the packing spectrum is proving the lower bound in Theorem \ref{Packing Main}. Unlike the Hausdorff dimension, the packing dimension of a level set typically exceeds the supremum of the dimensions of the invariant measures supported on that set. 
We construct a non-invariant measure specifically suited to obtaining an optimal lower bound for packing dimension.

The rest of the paper is structured as follows. We begin by restating Theorems \ref{Packing level sets} and \ref{Packing Main} in Section $2$ in terms of the symbolic space. The proof of Theorem \ref{Packing Main symbolic} is given in sections $3$ and $4$. In Section $3$ we prove the lower bound, and in Section $4$ we prove the upper bound. In Section $5$ we deduce some regularity properties of the packing spectrum. In Section $6$ we present two simple examples exhibiting some interesting features of the packing spectrum in the two dimensional case. In Section $7$ we conclude with some extensions of Theorem \ref{Packing level sets} which follow from Theorem \ref{Packing Main} along with some open questions.

\section{Symbolic Dynamics}
We begin by restating our theorem in terms of the associated symbolic space. Let $\Sigma$ denote the symbolic space $\D^{\N}$ under the usual product topology. We let $\Pi:\Sigma\rightarrow \Lambda$ denote the natural projection given by 
\begin{equation}
\Pi:(\omega_{\nu})_{\nu\in\N}\mapsto \left(\sum_{\nu\in\N}\frac{i_j(\nu)}{a_j^{\nu}}\right)_{j=1}^d\text{   where   }\omega_{\nu}=(i_j(\nu))_{j=1}^d\text{  for  } \nu\in\N
.\end{equation}
For each $k=1,\cdots,d$ we let $\eta_k$ denote the projection of $\D$ by $\eta_k:(i_j)_{j=1}^d \mapsto (i_j)_{j=1}^{d-(k-1)}$ and $\Sigma_k:=\eta_k(\D)^{\N}$. We then define a projection $\chi_k:\Sigma \mapsto \Sigma_k$, corresponding to $\pi_k:\T^d\rightarrow \T^{d-(k-1)}$ by $\chi_k:(\omega_{\nu})_{\nu=1}^{\infty}\mapsto (\eta_k(\omega_{\nu}))_{\nu=1}^{\infty}$. We let $\Pi_k:\Sigma_k\rightarrow \pi_k(\Lambda)$ denote the natural projection given by 
\begin{equation}
\Pi_k:(\tau_{\nu})_{\nu\in\N}\mapsto \left(\sum_{\nu\in\N}\frac{i_j(\nu)}{a_j^{\nu}}\right)_{j=1}^{d-(k-1)}\text{   where   }\tau_{\nu}=(i_j(\nu))_{j=1}^{d-(k-1)}\text{  for  } \nu\in\N
.\end{equation}
Note that $\Pi_k\circ\chi_k=\pi_k\circ\Pi$. We let $\sigma$ denote the left shift on $\Sigma$ and for each $k$, $\sigma_k$ denotes the left shift on $\Sigma_k$. Note that $f\circ\Pi=\Pi \circ \sigma$ and for each $k$ $f_k\circ\Pi_k=\Pi_k\circ \sigma_k$. 
Given a finite sequence $(\omega_{\nu})_{\nu=1}^n\in\eta_k(\D)$ we let $[\omega_1\cdots\omega_n]$ denote the cylinder set \begin{equation}
[\omega_1\cdots\omega_n]:=\left\lbrace \omega'\in\Sigma_k:\omega'_{\nu}=\omega_{\nu}\text{ for }\nu=1,\cdots,n\right\rbrace.
\end{equation} 
 Given $\varphi:\Sigma\rightarrow \R$ and $n\in \N$ we define $\var_n(\varphi)$ by
\begin{equation}
\var_n(\varphi):=\sup\left\lbrace |\varphi(\omega)-\varphi(\tau)|: \omega_{\nu}=\tau_{\nu} \text{ for }\nu=1,\cdots,n\right\rbrace
.\end{equation}
We also define $A_n(\varphi):\Sigma\rightarrow \R$ to be the map $\omega\mapsto \frac{1}{n}\sum_{l=0}^{n-1}\varphi(\sigma^l\omega)$.

We are interested in the space of all Borel probability measures $\M(\Sigma)$ under the weak $*$ topology. Since $\Sigma$ is compact and hence the space $C(\Sigma)$ of continuous real valued functions on $\Sigma$ is separable, we may choose a countable family of potentials $(\varphi_l)_{l \in \N}$ with norm one, $||\varphi_l||_{\infty}=1$, for all $l\in \N$, for which sets of the form 
\begin{equation}
\left\lbrace \nu \in \M(\Sigma): \bigg| \int \varphi_l d\nu - \int \varphi_l d\mu \bigg|< \epsilon \text{ for all }l\leq L \right\rbrace,
\end{equation}
with $\mu \in \M(\Sigma)$ and $L \in \N$, form a neighbourhood basis of $\M(\Sigma)$.

For each $n\in \N$ we let $\M_{\sigma^n}(\Sigma)$ denote the set of $\sigma^n$-invariant Borel probability measures, let $\E_{\sigma^n}(\Sigma)$ denote the set of $\mu\in \M_{\sigma^n}(\Sigma)$ which are ergodic, with respect to $\sigma^n$, and let $\B_{\sigma^n}(\Sigma)$ denote the set of $\mu\in \E_{\sigma^n}(\Sigma)$ which are also Bernoulli.

Given a Borel probability measure $\mu \in \M(\Sigma)$ we define 
\begin{equation}
A_n(\mu):= \frac{1}{n}\sum_{l=0}^{n-1} \mu \circ \sigma^{-l}.
\end{equation}
Given $\omega \in \Sigma$ we let $\V(\omega)$ denote the set of all weak $*$ accumulation points of the sequence of measures $(A_n(\delta_{\omega}))_{n\in\N}$ where $\delta_{\omega}$ denotes the Dirac measure concentrated at $\omega$. Given $A \subseteq \M_{\sigma}(\Sigma)$ we define 
\begin{eqnarray}
\Gamma(A):&=& \left\lbrace \omega \in \Omega : \V(\omega)=A\right\rbrace \\
\Omega(A):&=& \left\lbrace \omega \in \Omega : \V(\omega)\subseteq A\right\rbrace \nonumber.
\end{eqnarray}
For each $k \leq d$ we define $H^k(\sigma,A)$ by  
\begin{eqnarray*}
H^k(\sigma,A)&:=&\sup \left\lbrace h_{\mu\circ \chi_k^{-1}}(\sigma_k): \mu \in A \right\rbrace.
\end{eqnarray*}
We shall prove the following Theorem which implies Theorems \ref{Packing level sets} and \ref{Packing Main}.
\begin{theorem}\label{Packing Main symbolic}
Suppose that $A$ is a non-empty closed convex subset of $\M_{\sigma}(\Sigma)$. Then 
\begin{equation*}
\Dim \Pi(\Gamma(A))=\Dim \Pi(\Omega(A))=\frac{H^1(\sigma,A)}{\log a_1} +\sum_{k=2}^d \left(\frac{1}{\log a_k}-\frac{1}{\log a_{k-1}}\right) H^k(\sigma,A).
\end{equation*}
\end{theorem}

\section{Proof of the lower estimate}
Fix a non-empty closed convex subset $A\subseteq \M(\Sigma,\sigma)$. Take $\zeta>0$ and choose some $\mu_j\in\M_{\sigma}(\Sigma)$ for $j=1,\cdots, d$ such that
\begin{eqnarray} \label{mu_j def}
h_{\mu_j\circ \chi_j^{-1}}(\sigma_j)>H^j(\sigma,A) -\zeta.
\end{eqnarray}
Through a series of lemmas we shall prove that 
\begin{equation}
\Dim\Pi(\Gamma(A))\geq  \frac{H^1(\sigma,A)}{\log a_1} +\sum_{j=2}^d \left(\frac{1}{\log a_j}-\frac{1}{\log a_{j-1}}\right) H^j(\sigma,A).\end{equation}
To this end we construct a measure allowing us to apply the following result from geometric measure theory.
\begin{prop} \label{Measure Dimension Lemma}
Let $E\subseteq \R^n$ be a Borel set and $\mu$ a finite Borel measure. If $\limsup_{r\rightarrow 0}\frac{\log \mu(B(x;r))}{\log r}\geq s$ for all $x\in E$ and $\mu(E)>0$ then $\Dim(E)\geq s$.
\end{prop}
\begin{proof}
This follows from \cite{Falconer} Proposition 4.9. 
\end{proof}
Let $\lambda_0:=0$ and for $j=1,\cdots,d$ we let $\lambda_j:=\log a_d/\log a_j$. In order to obtain an optimal lower bound we shall construct a measure $\W$ which, for infinitly many values of $n$, behaves like $\mu_j$ for the digits from $\lceil \lambda_{j-1} n \rceil+1$ up to $\lceil \lambda_{j} n \rceil$, for each $j=1,\cdots, d$, and use this property to show that $\W \circ \Pi^{-1}$ has the required packing dimension. 

We must also choose $\W$ so that $\V(\omega)=A$ on a set of large $\W$ measure. To do this we take a sequence of measures $(m_q)_{q \in 2\N}$ in $A$ for which the set of weak $*$ limit points is of $(m_q)_{q \in 2\N}$ is precisely the set $A$. We shall also construct $\W$ so that, along a subsequence of times, $\W$ behaves like $(m_q)_{q\in 2\N}$.
 
To obtain such a measure, $\W$, we effectively piece together the various invariant measures that $\W$ is required to imitate. 
In order to carry out this procedure we must first approximate each of our invariant measures by members of $\bigcup_{n\in\N}\B_{\sigma^n}(\Sigma)$. This allows us to deal with three issues. Firstly, the invariant measures which $\W$ is required to mimic need not be ergodic. Nonetheless, there approximations will be ergodic for some $n$-shift $\sigma^n$, and this allows us to apply both Birkhoff's ergodic Theorem and the Shannon-McMillan-Breiman Theorem. Secondly, we do not assume King's disjointness condition (see \cite{King}) and allow our approximate squares to touch at their boundaries. As such we must insure that our measure is not too concentrated so that it behaves well under projection by $\Pi$. For members of $\bigcup_{n\in\N}\B_{\sigma^n}(\Sigma)$ we may do this simply by tweaking our measure so that it gives each finite word some positive probability. Thirdly, the process of pieceing together measures is greatly simplified by only working with members of $\bigcup_{n\in\N}\B_{\sigma^n}(\Sigma)$. This approximation introduces an error, both in the  expected local entropy and expected Birkhoff averages. However, these error terms go to zero as the approximation improves, so by concatenating increasingly good approximations we will obtain a measure which not only behaves well at every given stage, but gives positive measure to the level set $\Gamma(A)$ and gives an optimal lower bound for the packing dimension of $\Pi(\Gamma(A))$.

Similar techniques appear in the work of Gelfert and Rams \cite{Gelfert Rams}, Barral and Feng \cite{Barral Feng}, Baek, Olsen and Snigreva \cite{Snigreva Packing} and Barreira and Schmeling \cite{Nontypical1}.

\begin{lemma} \label{q Bernoulli convergence} For each $j=1,\cdots, d$ and $q\in2\N+1$ we may find $k(q)\in\N$ and $\nu^q_j\in\B_{\sigma^{k(q)}}(\Sigma)$ such that for $l=1\cdots, q$
\begin{enumerate}
\item [(i)] $ \displaystyle \big|h(\mu_j^q\circ \chi_j^{-1},\sigma_j)-\frac{1}{k(q)}h(\nu^q_j\circ \chi_j^{-1},\sigma_j^{k(q)})\big|<\frac{1}{q}$,
\vspace{4mm}
\item[(ii)] $ \displaystyle \big| \int A_{k(q)}(\varphi_l) d\nu^q_j- \int \varphi_l d\nu_j \big|<\frac{1}{q}$,
\vspace{4mm}
\item [(iii)] $ \displaystyle \nu^q_j([\omega_1\cdots\omega_{k(q)}])>0$ for all $(\omega_1,\cdots,\omega_{k(q)})\in\D^{k(q)}.$ 
\end{enumerate}
\end{lemma}
\begin{proof}
Given $j \in \N$ and $k\in\N$ we let $\mu^k_j$ denote the unique member of $\B_{\sigma^{k}}(\Sigma)$ which agrees with $\mu_j$ on cylinders of length $k$. So for all $(\tau_1,\cdots,\tau_k)\in\eta_j(\D)^k$ we have 
\begin{eqnarray*}
\mu^k_j\circ \chi_j^{-1}([\tau_1\cdots\tau_k])&=&\mu_a\circ\chi_j^{-1}([\tau_1\cdots\tau_k]).
\end{eqnarray*}
Now by the Kolmogorov-Sinai Theorem (\cite[Theorem 4.18]{Walters} ) we have 
\begin{eqnarray*}
 h(\mu_j\circ \chi_j^{-1},\sigma_j)=-\lim_{k\rightarrow \infty}\frac{1}{k}\sum_{(\tau_1,\cdots,\tau_k)\in\eta(\D)^k} \mu_j\circ \chi_j^{-1}([\omega_1,\cdots,\omega_k])\log \mu_j\circ \chi_j^{-1}([\omega_1,\cdots,\omega_k]).
\end{eqnarray*}  
Equivalently, 
\begin{eqnarray*}
 h(\mu_j \circ \chi_j^{-1},\sigma_j)&=&\lim_{k\rightarrow \infty}\frac{1}{k} h(\mu^k_j\circ \chi_j^{-1},\sigma^k).
 \end{eqnarray*}
Since each $\mu_j$ is $\sigma$ invariant we have $\int A_k(\varphi_l)d\mu_j=\int \varphi_l d\mu_j$ for $l=1,\cdots, k$ and as $\mu_j$ and $\mu^k_j$ agree on cylinders of length $k$ we have 
\begin{eqnarray*}
\bigg|\int A_k (\varphi_l) d\mu^k_j-\int A_k (\varphi_l)d\mu_j\bigg|&\leq& \frac{1}{k}
\sum_{n=0}^{k-1} \var_n (\varphi_l).
\end{eqnarray*}
Moreover, since each $\varphi_l$ is continuous $\var_k(\varphi_l)\rightarrow 0$ (and hence $\frac{1}{k}
\sum_{n=0}^{k-1} \var_n(\varphi_l)\rightarrow 0$) as $k\rightarrow\infty$. Thus, for each $l=1,\cdots,q$,
\begin{eqnarray*}
\lim_{k\rightarrow\infty} \int A_k(\varphi_l)d\mu_j^k=\int \varphi_l d \mu_j.
\end{eqnarray*}
Thus, taking $\nu^q_j=\mu^{k(q)}_j$, for each $j$, for sufficiently large $k(q)\in\N$ gives (i) and (ii). By slightly adjusting $\nu^q_j$ we may insure (i),(ii) and (iii) hold.
\end{proof}

\begin{lemma} \label{q Bernoulli convergence m} For all $q\in 2\N$ we may find $k(q)\in\N$ and $\tilde{m}_q\in\B_{\sigma^{k(q)}}(\Sigma)$ such that for $l=1\cdots, q$
\begin{enumerate}
\item[(i)] $ \displaystyle \big| \int A_{k(q)}(\varphi_l) d \tilde{m}_q- \int \varphi_l d m_q \big|<\frac{1}{q}$;
\vspace{4mm}
\item [(ii)] $ \displaystyle \tilde{m}_q([\omega_1\cdots\omega_{k(q)}])>0$ for all $(\omega_1,\cdots,\omega_{k(q)})\in\D^{k(q)}.$ 
\end{enumerate}
\end{lemma} 
\begin{proof} Essentially the same as Lemma \ref{q Bernoulli convergence}.
\end{proof}

Choose $\delta_{q}>0$ for each $q\in\N$ so that $\prod_{q\in\N}(1-\delta_{q})>0$. 

\begin{lemma} \label{N(q) convergence} For each $j=1,\cdots,d$ and $q\in 2\N+1$ we may find $N(q)\in\N$ and a subset $S^q_j\subseteq \Sigma$ with $\nu^q_j(S^q_j)>1-\delta_{q}$ and such that for all $\omega \in S^q_j$ and $n\geq N(q)$ and all $l=1,\cdots,q$ we have
\vspace{4mm}
\begin{enumerate}
\item [(i)] $ \displaystyle \bigg|\frac{1}{nk(q)}\sum_{r=0}^{nk(q)-1}\varphi_l(\sigma^r\omega)-\int \varphi_ld\mu_j \bigg|<\frac{1}{q}$,
\vspace{4mm}
\item [(ii)] $ \displaystyle \bigg|\frac{1}{nk(q)}\log\nu^q_j\circ \chi_j^{-1}([\eta_j(\omega_1)\cdots \eta_j(\omega_{nk(q)})])+h(\mu_j \circ \chi_j^{-1},\sigma_j) \bigg|<\frac{1}{q}$,
\vspace{4mm}
\item [(iii)] $ \displaystyle \left\lbrace d\in\D: \omega_r=d \text{ for some }r\leq N(q)\right\rbrace=\D$.
 \end{enumerate}
\vspace{4mm}
\end{lemma}

\begin{proof}
Given $q\in\N$ we may apply the Birkhoff ergodic theorem and the the Shannon-Breiman-MacMillan theorem to $\nu_j^q\circ \chi_j^{-1}\in\B_{\sigma_j^{k(q)}}(\Sigma_j)\subseteq \E_{\sigma_j^{k(q)}}(\Sigma_j)$ to obtain
\begin{equation}
\label{Birkhoff Con}
\lim_{n\rightarrow \infty} \frac{1}{nk(q)}\sum_{r=0}^{nk(q)-1} \varphi_l(\sigma^r \omega)=\lim_{n\rightarrow \infty} \frac{1}{n}\sum_{r=0}^{n-1}A_{k(q)}(\varphi_l)(\sigma^{rk(q)}\omega)=\int \varphi_l d\nu^q_j
\end{equation}
\begin{equation}\label{SMB Con}
\lim_{n\rightarrow \infty}\frac{1}{n}\log\mu^q_j \circ \chi_j^{-1}([\eta_j(\omega_1)\cdots \eta(\omega_{nk(q)})])=-h(\mu_j^q\circ \chi_j^{-1},\sigma_j^{k(q)})
\end{equation}
for $\nu_j^q$ almost every $\omega\in\Sigma_j$.

By Egorov's theorem we may choose subsets $S_j^q\subseteq \Sigma$ with $\nu^q_j(S^q_j)>1-\delta_{q}$ so that the convergences in (\ref{Birkhoff Con}) and (\ref{SMB Con}) are uniform on $S^q_j$. Thus, by Lemma \ref{q Bernoulli convergence} we choose $N(q)\in\N$ so that for all $n\geq N(q)$ and all $\omega\in S^q_j$ we have
\begin{eqnarray}
\label{Birkhoff a}
 &\bigg|\frac{1}{nk(q)}&\sum_{r=0}^{nk(q)-1}\varphi_l(\sigma^r\omega)-\alpha \bigg|<\frac{1}{q}\\
\label{SMB a}
 &\bigg|\frac{1}{nk(q)}&\log\nu^q_a\circ \chi_j^{-1}([\eta(\omega_1)\cdots\eta(\omega_{nk(q)})])+h(\mu_j\circ \chi_j^{-1},\sigma_j) \bigg|<\frac{1}{q}.
\end{eqnarray}
In light of condition (iii), for $\nu^q_j$ almost every $\omega\in\Sigma$ we have
$\{ d\in\D: \omega_r=d \text{ for some }r\in\N\}=\D$. Equivalently, for $\nu^q_j$ almost every $\omega\in\Sigma$ there exists some $M(\omega)\in\N$ for which
$\left\lbrace d\in\D: \omega_r=d \text{ for some }r\leq M(\omega)\right\rbrace=\D$. Thus, by moving to subset of $S^q_j$ of large $\nu^q_j$ measure, and increasing $N(q)$, if necessary, we may assume that for all $\omega\in S^q_j$ and all $n\geq N(q)$ we have 
\begin{equation}\label{every digit}
\left\lbrace d\in\D: \omega_l=d \text{ for some }l\leq N(q)\right\rbrace=\D.
\end{equation}
\end{proof}

\begin{lemma}\label{N(q) convergence 2} For each $q\in 2\N$ we may find $N(q)\in\N$ and a subset $S^q\subseteq \Sigma$ with $\tilde{m}^q(S^q)>1-\delta_{q}$ and such that for all $\omega \in S^q$ and $n\geq N(q)$ and all $l=1,\cdots,q$ we have
\vspace{4mm}
\begin{enumerate}
\item [(i)] $ \displaystyle \bigg|\frac{1}{nk(q)}\sum_{r=0}^{nk(q)-1}\varphi_l(\sigma^r\omega)-\int \varphi_ld m_q \bigg|<\frac{1}{q}$,
\vspace{4mm}
\item [(ii)] $ \displaystyle \left\lbrace d\in\D: \omega_r=d \text{ for some }r\leq N(q)\right\rbrace=\D$.
 \end{enumerate}
\vspace{4mm}
\end{lemma}
\begin{proof} Essentially the same as Lemma \ref{N(q) convergence}.
\end{proof}

We shall now construct a probability measure $\W$ on $\Sigma$. To do this we first define a rapidly increasing sequence of natural numbers $(\gamma_q)_{q\in\N}$ as follows. Let $\gamma_0:=0$ and for each $q\geq 1$ taking some $\gamma_q> (q+1)\prod_{j=1}^d (\lambda_j-\lambda_{j-1})^{-1}  (\prod_{r=1}^{q+2}N(r)k(r) +\gamma_{q-1})$ so that $\gamma_q-\gamma_{q-1}$ is divisible by $k(q)$. For each $k=1,\cdots, d$ we sequences of natural numbers $(\vartheta^k_q)_{q\in 2\N+1}$ by letting $\vartheta^k_q$ denote the greatest integer which is divisible by $k(q)$ and does not exceed $\lambda_k \gamma_q$. For simplicity we also let $\vartheta^{0}_{q}:=\gamma_{q-1}$. 
 
We define a measure $\W$ on $\Sigma$ by first defining $\W$ on cylinders of length $\gamma_{2Q}$ for some $Q\in \N$ and then extending $\W$ to a Borel probability measure via the Daniell-Kolmogorov consistency theorem (see \cite[Section 0.5]{Walters} ). Given a cylinder  $[\omega_1\cdots\omega_{\gamma_{2Q}}]$ of length $\gamma_{2Q}$ we let
\begin{equation}
\W([\omega_1\cdots\omega_{\gamma_{2Q}}]):=\hspace{8cm} 
\end{equation}
\begin{eqnarray}
\hspace{2cm}
\prod_{q=1}^{Q}\bigg(\prod_{j=1}^d\nu^{2q-1}_j ([\omega_{\vartheta^{j-1}_{2q-1}+1} \cdots \omega_{\vartheta^{j}_{2q-1}}]) \times \tilde{m}^{2q}([\omega_{\gamma_{2q-1}+1}\cdots, \omega_{\gamma_{2q}}])\bigg)
\nonumber.\end{eqnarray}
Define $S\subseteq \Sigma$ by,
\begin{eqnarray}
S:=\bigcap_{q=1}^{\infty}\bigg(\bigcap_{j=1}^d \left\lbrace \omega \in \Sigma: [\omega_{\vartheta^{j-1}_{2q-1}+1} \cdots \omega_{\vartheta^{j}_{2q-1}}] \cap S^{2q-1}_j \neq \emptyset \right\rbrace \\
\cap \left\lbrace \omega \in \Sigma : [\omega_{\gamma_{2q-1}+1}\cdots, \omega_{\gamma_{2q}}]\cap S^{2q}\neq \emptyset\right\rbrace\bigg).\nonumber
\end{eqnarray}
\begin{lemma} \label{S has positive W measure}$\W(S)>0$.
\end{lemma}
\begin{proof}
\begin{equation}\W(S)\geq \prod_{q=1}^{\infty}\left(\left(\prod_{j=1}^d \nu_j^{2q-1}(S^{2q-1}_j)\right)\tilde{m}^{2q}(S^{2q})\right)>\prod_{q=1}^{\infty}(1-\delta_q)^d>0. \end{equation}
\end{proof}

\begin{lemma}\label{V containment} For all $\omega \in S$, $\V(\omega)\subseteq A$.
\end{lemma}
\begin{proof} Choose $\omega \in S$ and fix $Q \in \N$ and $\epsilon>0$. For each $q \in 2\N$ with $q \geq Q$ take $\tau^q \in [\omega_{\gamma_{q-1}+1}\cdots, \omega_{\gamma_{q}}]\cap S^{q}$, which is non-empty since $\omega \in S$. By Lemma \ref{N(q) convergence 2}, for all $n\geq k(q)N(q)$,
\begin{equation}
\bigg|\sum_{r=0}^{n-1}\varphi_l(\sigma^r\tau^q)-n\int \varphi_l d m_q \bigg|<\frac{n}{q}+k(q).
\end{equation}
Hence, for all $N(q)k(q)\leq n \leq \gamma_q-\gamma_{q-1}$,
\begin{eqnarray}\label{ est1 }
\bigg|\sum_{r=\gamma_{q-1}}^{\gamma_{q-1}+n-1}\varphi_l(\sigma^r\omega)-n\int \varphi_l d m_q \bigg|<\sum_{r=0}^n\var_r(\varphi_l)+\frac{n}{q}+k(q).
\end{eqnarray}
In a similar way we can show that for all $q \in 2\N-1$, with $q\geq Q$, $j=1,\cdots, d$ and all $N(q)k(q)\leq n \leq \vartheta^{j-1}_q-\vartheta^j_{q-1}$ we have
\begin{eqnarray}\label{ est2 }
\bigg|\sum_{r=\vartheta^{j-1}_q}^{\vartheta^{j-1}_q+n-1}\varphi_l(\sigma^r\omega)-n\int \varphi_l d \mu_j \bigg|<\sum_{r=0}^n\var_r(\varphi_l)+\frac{n}{q}+k(q).
\end{eqnarray}
Moreover, given any $n,k,q \in \N$ we automatically have 
\begin{eqnarray}\label{ est3 }
\bigg|\sum_{r=k}^{k+n}\varphi_l(\sigma^r\omega)-n\int \varphi_l d m_q \bigg|<n.
\end{eqnarray}
Suppose $\gamma_{2q}<N <\gamma_{2q+2}$ where $2q-2 \geq Q$. Now consider the sum $\sum_{r=0}^{N-1}\varphi_l(\sigma^r\omega)$, for $l\leq Q$. First break the sum down as follows,
\begin{equation}
\sum_{r=0}^{N-1}\varphi_l(\sigma^r\omega)= \sum_{r=0}^{\gamma_{2q}-1}\varphi_l(\sigma^r\omega) +\sum_{r=\gamma_{2q}}^{N-1}\varphi_l(\sigma^r\omega). 
\end{equation}
To deal with the first summand, $\sum_{r=0}^{\gamma_{2q}-1}\varphi_l(\sigma^r\omega)$, we write,
\begin{equation}\label{ est11 }
\sum_{r=0}^{\gamma_{2q}-1}\varphi_l(\sigma^r\omega)=\underbrace{\sum_{r=0}^{\gamma_{2q-2}-1}\varphi_l(\sigma^r\omega)}_{*}+ \sum_{j=1}^d \underbrace{\sum_{r=\vartheta^{j-1}_{2q-1}}^{\vartheta^{j}_{2q-1}-1}\varphi_l(\sigma^r\omega)}_{**} +\underbrace{\sum_{r=\gamma_{2q-1}}^{\gamma_{2q}-1}\varphi_l(\sigma^r\omega)}_{***}.
\end{equation}
To part $*$ we apply (\ref{ est3 }) whilst to each of the parts labeled $**$ we apply (\ref{ est2 }) and to the part labeled $***$ we apply (\ref{ est1 }).
 
For the second summand, $\sum_{r=\gamma_{2q}}^{N-1}\varphi_l(\sigma^r\omega)$, there are two cases. Either we have $N \leq \gamma_{2q+1}$ or $N>\gamma_{2q+1}$. In the former case we have
\begin{equation}\label{ est22 }
\sum_{r=\gamma_{2q}}^{N-1}\varphi_l(\sigma^r\omega)=\sum_{j=1}^{J}\underbrace{\sum_{r=\vartheta^{j-1}_{2q}}^{\vartheta^{j}_{2q}-1}\varphi_l(\sigma^r\omega)}_{**}+ \underbrace{\sum_{r=\vartheta^{J}_{2q}}^{N-1}\varphi_l(\sigma^r\omega)}_{\dagger}, 
\end{equation}
where $J$ is the greatest $j \in \{1, \cdots, d\}$ such that $J<N-1$. To parts labeled $**$ we again apply (\ref{ est2 }), and to the part labeled ($\dagger$) we either apply (\ref{ est3 }) or (\ref{ est2 }), depending on whether 
$N- \vartheta^{J}_{2q}-1<  \max_{j \in \{2q,2q+1\}}\{k(j)N(j)\}$ or $N-\vartheta^{J}_{2q}-1\geq  \max_{j \in \{2q,2q+1\}}\{k(j)N(j)\}$. In the latter case we have,
\begin{equation}\label{ est33 }
\sum_{r=\gamma_{2q}}^{N-1}\varphi_l(\sigma^r\omega)=\sum_{j=1}^{d}\underbrace{\sum_{r=\vartheta^{j-1}_{2q}}^{\vartheta^{j}_{2q}-1}\varphi_l(\sigma^r\omega)}_{**}+ \underbrace{\sum_{r=\gamma_{2q+1}}^{N-1}\varphi_l(\sigma^r\omega)}_{\dagger\dagger}. 
\end{equation}
Again, to the parts labeled $**$ we apply (\ref{ est2 }), and to the part labeled ($\dagger\dagger$) we either apply (\ref{ est3 }) or (\ref{ est1 }), depending on whether 
$N- \gamma_{2q+1}-1<  \max_{j \in \{2q,2q+1\}}\{k(j)N(j)\}$ or $N-\gamma_{2q+1}-1\geq  \max_{j \in \{2q,2q+1\}}\{k(j)N(j)\}$.
 
Thus, by combining (\ref{ est11 }), (\ref{ est22 }), (\ref{ est33 }), in each case we see that there exists $\beta_{2q}^N, \beta_{2q+2}^N,\lambda_{1}^N, \cdots, \lambda_{d}^N  \in [0,1]$ which sum to one, depending solely on $N$ and not on $l=1,\cdots,Q$, for which we have 
\begin{eqnarray*}
\bigg|\sum_{r=0}^{N-1}\varphi_l(\sigma^r\omega)- \sum_{j\in \{2q, 2q+2\}}N \beta^N_j\int \varphi_l dm_j - \sum_{j=1}^d N \lambda^N_j  \int \varphi_l d \mu_j \bigg|\\
<(2d+1)\left(\sum_{r=0}^N\var_r(\varphi_l)+\frac{N}{q-1}+\max_{j \in \{0,1,2,3\}}\{k(2q+j)\}\right)\\+\left(\gamma_{2q-2}+\max_{j \in \{2q,2q+1\}}\{k(j)N(j)\}\right)\\
< o(N)+o(\gamma_{2q})\leq o(N),
\end{eqnarray*}
where we use the continuity of each $\varphi_l$ together with the definition of $(\gamma_q)_{q\in\N}$ to obtain the last line. Moreover, since $A$ is convex, for each such $N$ the measure $\rho_N:=\sum_{j\in \{2q, 2q+2\}}N \beta^N_j m_j + \sum_{j=1}^d \lambda^N_j \mu_j$ is a member of $A$. Hence, for all sufficiently large $N$, we have 
\begin{equation}
\bigg|\int \varphi_l dA_N(\delta_{\omega})- \int \varphi_l d\rho_N \bigg|<\epsilon,
\end{equation}
for $l=1,\cdots,Q$. Since $A$ is also closed it follows that every weak $*$ accumulation point of the sequence $(A_N(\delta_{\omega}))_{N\in\N}$ is a member of $A$.
\end{proof}
\begin{lemma}\label{A containment} For all $\omega \in S$, $A \subseteq \V(\omega)$.
\end{lemma}

\begin{proof} 
Take $\omega \in S$, $\alpha \in A$. Since the set of accumulation points of $(m_q)_{q\in 2\N}$ is equal to $A$ we may extract a subsequence $(m_{q_j})_{j\in\N}$ converging to $\rho$. Now choose $\epsilon>0$ and choose $Q$ so large that for $j$ with $q_j\geq Q$
\begin{eqnarray}
\bigg|\int\varphi_l dm_{q_j}-\int \varphi_l d \rho \bigg|<\epsilon.
\end{eqnarray}
As in \ref{ est1 } we see that for all $j$ with $q_j\geq Q$, $l=1,\cdots,Q$,
\begin{eqnarray}
\bigg|\sum_{r=\gamma_{q_j-1}}^{\gamma_{q_j}-1}\varphi_l(\sigma^r\omega)-({\gamma_{q_j}-\gamma_{q_j-1}})\int \varphi_l d m_{q_j} \bigg|\\<\sum_{r=0}^{\gamma_{q_j}-\gamma_{{q_j}-1}}\var_r(\varphi_l)+\frac{{\gamma_{q_j}-\gamma_{{q_j}-1}}}{{q_j}}+k({q_j}).
\end{eqnarray} 
Hence, for all ${q_j}\geq Q$, $l=1,\cdots,Q$,
\begin{eqnarray}
\bigg|\sum_{r=0}^{\gamma_{{q_j}}-1}\varphi_l(\sigma^r\omega)-\gamma_{q_j}\int \varphi_l d m_{q_j} \bigg|\\<\gamma_{{q_j}-1}+\sum_{r=0}^{\gamma_{q_j}-\gamma_{{q_j}-1}}\var_r(\varphi_l)+\frac{{\gamma_{q_j}-\gamma_{{q_j}-1}}}{{q_j}}+k({q_j}).
\end{eqnarray}
Thus, by the definition of $(\gamma_q)_{q\in \N}$ and the fact that each $\varphi_l$ is continuous, we have
\begin{eqnarray}
\bigg|\int\varphi_l dA_{\gamma_{q_j}}(\delta_{\omega})-\int \varphi_l d \rho \bigg|<2\epsilon.
\end{eqnarray}
It follows that $A_{\gamma_{q_j}}(\delta_{\omega})\rightarrow \rho$ as $j\rightarrow \infty$.
\end{proof}

\begin{lemma} \label{S is in Gamma} $S \subseteq \Gamma(A)$.
\end{lemma}
\begin{proof}
Combine Lemmas \ref{V containment} and \ref{A containment}.
\end{proof}

For each $q\in2\N+1$ we define the $q$th approximate square $B^L_q(\omega)$ to be the set
\begin{equation*}
B^L_q(\omega):=\left\lbrace \omega'\in \Sigma: \text{ For }j=1,\cdots,d \hspace{.2cm} \eta_j(\omega'_{\nu})=\eta_j(\omega_{\nu})\text{ for  }\nu=1,\cdots, \vartheta^j_q  \right\rbrace
.\end{equation*}
\begin{lemma} \label{Entropy lemma}
For all $\omega\in S$ 
\begin{equation*} 
\limsup_{q\rightarrow \infty}\frac{1}{\gamma_q}\log \W(B^L_q(\omega))\leq -\lambda_1 h(\mu_1,\sigma)-\sum_{j=1}^d(\lambda_j-\lambda_{j-1})h(\mu_j\circ \chi_j^{-1},\sigma_j)
.\end{equation*}
\end{lemma}
\begin{proof}
Take $\omega\in S$ and $q\in2\N+1$ and $j\in \{1, \cdots, d\}$. Since $\vartheta^j_{q}-\vartheta^{j-1}_{q}\geq N(q)$  and $[\omega_{\gamma_{q-1}+1}\cdots\omega_{\vartheta_q}]\cap S^q_j\neq \emptyset$ it follows from the definition of $\W$ and $(\gamma_q)_{q\in\N}$ together with Lemma \ref{N(q) convergence} that
 \begin{equation}
 \log\W\left(\sigma^{-\vartheta^{j-1}_{q}}[\eta_j(\omega_{\vartheta^{j-1}_{q}+1})\cdots\eta_j(\omega_{\vartheta^j_q})]\right)\hspace{3cm}
 \end{equation}
 \begin{eqnarray*}
 &=&\log\nu_j^q([\eta_j(\omega_{\vartheta^{j-1}_{q}+1})\cdots\eta_j(\omega_{\vartheta^j_q})])\\
&\leq& -\left(\vartheta^j_{q}-\vartheta^{j-1}_{q}\right)h(\mu_j^q\circ \chi_j^{-1},\sigma_j) +\frac{\vartheta^j_{q}-\vartheta^{j-1}_{q}}{q}
\\
&\leq& -\left((\lambda_j-\lambda_{j-1})\gamma_{q}-k(q)-\gamma_{q-1}\right)h(\mu_j^q\circ \chi_j^{-1},\sigma_j) +\frac{\gamma_q}{q}\\
&\leq& -(\lambda_j-\lambda_{j-1})\gamma_{q}h(\mu_j^q\circ \chi_j^{-1},\sigma_j) +o(\gamma_q).
\end{eqnarray*}
By the definition of $\W$, each of the cylinders $\sigma^{-\vartheta^{j-1}_{q}}[\eta_j(\omega_{\vartheta^{j-1}_{q}+1})\cdots\eta_j(\omega_{\vartheta^j_q})]$, with $j=1,\cdots,d$, are independent with respect to $\W$. Thus, letting $q\rightarrow \infty$ we have
\begin{eqnarray} 
\limsup_{q\rightarrow \infty}\frac{1}{\gamma_q}\log \W\left(\bigcap_{j=1}^d \sigma^{-\vartheta^{j-1}_{q}}[\eta_j(\omega_{\vartheta^{j-1}_{q}+1})\cdots\eta_j(\omega_{\vartheta^j_q})] \right)\\
\leq -\lambda_1 h(\mu_1,\sigma)-\sum_{j=1}^d(\lambda_j-\lambda_{j-1})h(\mu_j\circ \chi_j^{-1},\sigma_j)\nonumber.
\end{eqnarray}
Since $B^L_q(\omega)\subseteq\bigcap_{j=1}^d \sigma^{-\vartheta^{j-1}_{q}}[\eta_j(\omega_{\vartheta^{j-1}_{q}+1})\cdots\eta_j(\omega_{\vartheta^j_q})]$ the lemma follows.
\end{proof}

The following lemma allows us to deal with the fact that our approximate squares may meet at their boundaries.
\begin{lemma} \label{Containment}
For all $\omega\in S$ and $q\in\N$ we have
\[B\left(\Pi(\omega);\min_{j\in\{1,\cdots,d\}}\{a_j^{-\vartheta^j_q-\max\{N(q),N(q+1)\}}\}\right)\cap \Lambda \subseteq \Pi(B^L_q(\omega)).\]
\end{lemma}
\begin{proof}
Fix $\omega=((i^j_{\nu})_{j=1}^d)_{\nu\in\N}\in S$ and let $M(q):= \max\{N(q),N(q+1)\}$. Clearly it suffices to show that for each $j=1,\cdots, d$,
\begin{equation}\label{1 d ball}
B\left(\Pi_j(\omega);\min_{j\in\{1,\cdots,d\}}\{a_j^{-\vartheta^j_q-M(q)}\}\right)\cap \Lambda \subseteq \Pi_j(B^L_q(\omega)).
\end{equation}
Take $j\in\{1,\cdots,d\}$ and let $x=\Pi_j(\omega)$. 
We may divide $\Pi_j(B^L_q(\omega))$ into $\#\eta_j(\D)^{M(q)}$ intervals of width $a_j^{-\vartheta^j_q-M(q)}$ with disjoint ineriors, each corresponding to a possible string of digits $i'_{\vartheta^j_q+1}\cdots i'_{\vartheta^j_q+M(q)}$ for $\tau=((\tilde{i}^j_{\nu})_{j=1}^d)_{\nu\in\N}\in B^L_q(\omega)$. Since $\omega\in S$ we have $[\omega_{\vartheta^j_q+1}\cdots\omega_{\vartheta^j_q+N(q)}]\cap S^q_{j+1}\neq \emptyset$ for $j=1,\cdots, d-1$ and $[\omega_{\vartheta^d_q+1}\cdots\omega_{\vartheta^d_q+N(q+1)}]\cap S^{q+1}_d\neq \emptyset$. Thus, by Lemma \ref{N(q) convergence} in the first case and Lemma \ref{N(q) convergence 2}
\begin{eqnarray*}
\D&=&\left\lbrace d\in\D: \omega_l=d \text{ for some }\vartheta^j_q< l\leq \vartheta^j_q+M(q)\right\rbrace.
\end{eqnarray*}
It follows from $\eta_j(\D)=\left\lbrace d\in\eta_j(\D): i^j_{\nu}=d \text{ for some }\vartheta^j_q< l\leq \vartheta^j_q+N(q)\right\rbrace$ that $x$ is in neither the far left nor the far right interval of $\Pi_j(B^L_q(\omega))$, for in either case $(i^j_{\nu})_{\nu=\vartheta^j_q+1}^{\vartheta^j_q+N(q)}$ would be a constant sequence. Since $\#\eta_j(\D)>1$ it follows that $x$ is a distance at least $a_j^{-\vartheta^j_q-M(q)}$ from any point $y$ such that $y$ as an $a_j$-ary digit expansion $(\tilde{i}^j_{\nu})_{\nu\in\N}$ with $\tilde{i}^j_{\nu}\neq i^{j}_{\nu}$ for some $\nu \leq \vartheta^j_q$. Thus, (\ref{1 d ball}) holds.
\end{proof}
Let $\M:=\W\circ\chi^{-1}$ denote the pushdown of $\W$ onto $\Lambda$.

\begin{lemma}\label{local dimension lower bound} For all $\omega\in S$ 
\[\limsup_{r\rightarrow 0} \frac{\log \M(B(\chi(\omega);r))}{\log r}\geq -\frac{h(\mu_1,\sigma)}{\log a_1}+\sum_{j=1}^d\left(\frac{1}{\log a_j}-\frac{1}{\log a_{j-1}}\right)h(\mu_j\circ \chi_j^{-1},\sigma_j).\]
\end{lemma}
\begin{proof}
Recall that for each $j=1,\cdots, d$ we defined $\lambda_j:=\log a_d/\log a_j$, and for each $q\in 2\N+1$ we have $\vartheta_q^j<\lambda_j\gamma_q$ and so
\begin{equation}
\min_{j\in\{1,\cdots,d\}}\{a_j^{-\vartheta^j_q-\max\{N(q),N(q+1)\}}\}\geq a_d^{-\gamma_q-\max\{N(q),N(q+1)\}/\lambda_1}.
\end{equation}
Choose $\omega \in S$. By Lemma \ref{Containment}, 
\begin{equation}
B\left(\Pi(\omega); a_d^{-\gamma_q-\max\{N(q),N(q+1)\}/\lambda_1} \right)\cap \Lambda \subseteq \Pi(B^L_q(\omega)).
\end{equation}
Hence,
\begin{equation}
\M\left(B\left(\Pi(\omega); a_d^{-\gamma_q-\max\{N(q),N(q+1)\}/\lambda_1} \right)\right)\leq \W\left(B^L_q(\omega)\right).
\end{equation}
It follows from Lemma \ref{Entropy lemma} that
\begin{eqnarray} 
\limsup_{q\rightarrow \infty}\frac{1}{\gamma_q}\log \M\left(B\left(\Pi(\omega); a_d^{-\gamma_q-\max\{N(q),N(q+1)\}/\lambda_1} \right)\right)\\
\leq -\lambda_1 h(\mu_1,\sigma)-\sum_{j=2}^d(\lambda_j-\lambda_{j-1})h(\mu_j\circ \chi_j^{-1},\sigma_j)
\nonumber.\end{eqnarray}
Thus, noting that $\max\{N(q),N(q+1)\}=o(q)$, by the definition of $(\gamma_q)_{q\in \N}$ we have 
\begin{eqnarray}
\limsup_{q\rightarrow \infty} \frac{\log \M\left(B\left(\chi(\omega);a_d^{-\gamma_q-\max\{N(q),N(q+1)\}/\lambda_1}\right)\right)}{\log a_d^{-\gamma_q-\max\{N(q),N(q+1)\}/\lambda_1}}
\\ \geq -\frac{h(\mu_1,\sigma)}{\log a_1}+\sum_{j=2}^d\left(\frac{1}{\log a_j}-\frac{1}{\log a_{j-1}}\right)h(\mu_j\circ \chi_j^{-1},\sigma_j)\nonumber.
\end{eqnarray}

\end{proof}
Since $\M(\Pi(S))\geq \W(S)>0$ we may combine Proposition \ref{Measure Dimension Lemma} with Lemma \ref{local dimension lower bound} to see that
\begin{equation}
\Dim(\Pi(S))\geq -\frac{h(\mu_1,\sigma)}{\log a_1}+\sum_{j=1}^d\left(\frac{1}{\log a_j}-\frac{1}{\log a_{j-1}}\right)h(\mu_j\circ \chi_j^{-1},\sigma_j).
\end{equation}
Thus, by Lemma \ref{S is in Gamma} and our choice of $\mu_j$ (\ref{mu_j def}) we have 
\begin{equation}
\Dim\Pi(\Gamma(A))\geq  \frac{H^1(\sigma,A)}{\log a_1} +\sum_{j=2}^d \left(\frac{1}{\log a_j}-\frac{1}{\log a_{j-1}}\right) H^j(\sigma,A)- \zeta\frac{d}{\log a_1}.
\end{equation}
By letting $\zeta \rightarrow 0$ this completes the proof of the lower bound.
\section{Proof of the Upper Bound}
Take $A \subseteq \M_{\sigma}(\Sigma)$. Recall that for each $q \in \N$ we defined
\begin{equation}
U(A,q):=\left\lbrace \mu \in \M(\Sigma): \exists \nu \in A\hspace{2mm}\forall l \leq q \hspace{2mm} \bigg|\int \varphi_l d\mu-\int \varphi_l d\nu\bigg|\leq \frac{1}{q}\right\rbrace.
\end{equation}
Given $N\in\N$ we let
\begin{equation}
\Omega(A,N,q):= \left\lbrace \omega \in \Sigma :  \forall n \geq N \hspace{2mm} A_n(\delta_{\omega})\in U(A,q)\right\rbrace.
\end{equation}
Note that for each $q \in \N$ $\{\Pi(\Omega(A,N,q))\}_{N\in\N}$ is a countable cover of $\Pi(\Omega(A))$. As such we shall give an estimate for the upper box dimension of the sets $\Pi(\Omega(A,N,q))$ before applying the following reformulation of the notion of packing dimension.

\begin{prop} \label{BoxPacking}
Given $E\subseteq \R^n$ we have
\begin{equation*}
\Dim E=\inf\left\lbrace \sup_{n\in\N} \BDim E_n:E\subseteq \bigcup_{n\in\N}E_n\right\rbrace,
\end{equation*}
where the infimum is taken over all countable covers $\{E_n\}_{n\in\N}$ of $E$.
\end{prop}
The above formula is equivalent to the usual definition of packing dimension in terms of $s$-dimensional packing measures (see \cite[Section 5.9 and Theorem 5.11]{Mattila}).
 
Recall that for each $j=1,\cdots,d$, $\lambda_j:=\log a_d/\log a_j$. We also let $\lambda_0:=0$. Given $n\in \N$ we define $\A^j_n$ to be the set of all $(\tau_{\lceil \lambda_{j-1} n\rceil+1},\cdots,\tau_{\lceil \lambda_j n\rceil})\in \eta_j(\D)^{\lceil \lambda_j n\rceil-\lceil \lambda_{j-1} n\rceil}$ satisfying 
\begin{eqnarray} \chi_j\left(\Omega(A,N,q)\right) \cap \sigma_j^{-\lceil \lambda_{j-1} n\rceil}[\tau_{\lceil \lambda_{j-1} n\rceil+1},\cdots,\tau_{\lceil \lambda_j n\rceil}]\neq \emptyset.
\end{eqnarray}

\begin{lemma} \label{First Box Estimate}
\[\BDim{\Pi(\Omega(A,N,q))}\leq \sum_{j=1}^d  \limsup_{n\rightarrow \infty}\frac{\log \#\A^j_n}{n\log a_d}.\]
\end{lemma} 
\begin{proof}
Given $r>0$ we let $N(r)$ denote the minimal number of balls of radius $r$ required to cover $\Pi(\Omega(A,N,q))$ so that
\begin{equation*}
\BDim\Pi(\Omega(A,N,q))=\limsup_{r\rightarrow 0}\frac{\log N(r)}{\log(1/r)}
.\end{equation*}
For each $n_r$ we take $n_r\in \N$ so that $a_d^{-n_r}<r\leq a_d^{-n_r+1}$.
Given $\kappa:=(\kappa_j)_{j=1}^d$ where $\kappa_j=(\tau^j_{\lceil \lambda_{j-1} n\rceil+1},\cdots,\tau^j_{\lceil \lambda_j n\rceil})\in\A^j_{n_r}$ we let $B(\kappa)$ denote the approximate square
\begin{equation}
B(\kappa):=\Pi \left(\bigcap_{j=1}^d\sigma^{-\lceil \lambda_{j-1} n\rceil}\chi_j^{-1}\left([\tau_{\lceil \lambda_{j-1} n\rceil+1},\cdots,\tau_{\lceil \lambda_j n\rceil}]\right)\right)
.\end{equation}
It follows from the definition of $\Pi$ that each $B(\kappa)$ has diameter no greater than $a_j^{-\lceil \lambda_j n_r\rceil}\leq a_d^{-n_r}$. Moreover, 
\begin{equation}
\Pi(\Omega(A,N,q))\subseteq \bigcup_{\kappa\in\prod_{j=1}^d\A^j_n}B(\kappa)
.\end{equation}
Thus,
\begin{equation}
N(r)\leq \#\left\lbrace B(\kappa):\kappa\in\prod_{j=1}^d\A^j_{n_r}\right\rbrace =\prod_{j=1}^d \#\A^j_{n_r}.
\end{equation}
Hence, since $r\leq a_d^{-n_r+1}$,
\begin{eqnarray}
\limsup_{r\rightarrow 0}\frac{\log N(r)}{\log(1/r)}&\leq& \limsup_{r\rightarrow 0}\frac{\sum_{j=1}^d\log \#\A^j_{n_r}}{n_r\log a_d}\frac{ n_r}{n_r-1} \\
&\leq& \limsup_{n\rightarrow \infty}\frac{\sum_{j=1}^d\log \#\A^j_n}{n\log a_d} \nonumber \\
&\leq& \sum_{j=1}^d \limsup_{n\rightarrow \infty}\frac{\log \#\A^j_n}{n\log a_d}\nonumber.
\end{eqnarray}

\end{proof}

Recall that given $\nu \in \M_{\sigma^k}(\Sigma)$ for each $k\in\N$ we let $A_k(\nu):=1/k\sum_{l=0}^{k-1} \nu 
\circ \sigma^{-l}$.
\begin{lemma}\label{averagemeasure}
If $\mu=A_k(\nu)$ for some $\nu \in \E_{\sigma^k}(\Sigma)$ then for each $j=1,\cdots, d$ and all $\varphi\in C(\Sigma)$ 
\vspace{4mm}
\begin{enumerate}
\item [(i)] $\displaystyle\mu\circ \chi_{j}^{-1} \in\E_{\sigma_j}(\Sigma_j)$
\vspace{4mm}
\item [(ii)] $\displaystyle h(\mu \circ \chi_j^{-1},\sigma_j)=1/k h(\nu\circ \chi_j^{-1},\sigma_j^k)$
\vspace{4mm}
\item [(iii)] $\displaystyle \int \varphi d\mu=\int A_k(\varphi) d\nu$.
\end{enumerate}
\vspace{4mm}
\end{lemma}
\begin{proof}
By noting that $A_k(\nu)\circ \chi_j^{-1}=\sum_{r=1}^{k-1}\nu \circ \chi_j^{-1} \circ\sigma_j^{-r}$ we see that Lemma \ref{averagemeasure} follows lemma in \cite{non-uniformly hyperbolic} Lemma 2.
\end{proof}

For each $l\in\N$ we define   
\begin{eqnarray}
H^j(A,l)&:=&\sup \left\lbrace h(\mu\circ \chi_j^{-1},\sigma_j): \mu \in U(A,l) \right\rbrace.
\end{eqnarray}  
Define a constant $L\in(0,1)$ by
\begin{equation}
L:=\min_{j\in\{1,\cdots,d\}} \left\lbrace \frac{\lambda_j - \lambda_{j-1} }{4\lambda_j - \lambda_{j-1}}\right\rbrace.
\end{equation}
\begin{lemma} \label{Entropy Bound} For all $j=1,\cdots,d$,
\begin{eqnarray*}
\limsup_{n\rightarrow \infty}\frac{1}{n}\log \#\A^j_{n} &\leq& (\lambda_j-\lambda_{j-1}) H^j\left(A,\lceil Lq \rceil \right).\end{eqnarray*}
\end{lemma}
\begin{proof}
Take $j\in \{1,\cdots,d\}$. For each $\tau=(\tau_{\lceil \lambda_{j-1} n\rceil+1},\cdots,\tau_{\lceil \lambda_j n\rceil})\in\A^j_{n_r}$ choose $\kappa^{\tau}=(\omega^{\tau}_{\lceil \lambda_{j-1} n\rceil+1},\cdots,\omega^{\tau}_{\lceil \lambda_j n\rceil})\in \D^{\lceil \lambda_j n\rceil-\lceil \lambda_{j-1} n\rceil}$ so that 
\begin{eqnarray*}
\left(\eta_j(\omega^{\tau}_{\lceil \lambda_{j-1} n\rceil+1}),\cdots,\eta_j(\omega^{\tau}_{\lceil \lambda_{j} n\rceil})\right)=\tau,
\end{eqnarray*} 
and
\begin{eqnarray} \label{omega of tau def}
\Omega(A,N,q) \cap \sigma^{-\lceil \lambda_{j-1} n\rceil}[\omega^{\tau}_{\lceil \lambda_{j-1} n\rceil+1}\cdots \omega^{\tau}_{\lceil \lambda_{j} n\rceil}]\neq \emptyset.
\end{eqnarray}
We now let $\nu_n\in \B_{\sigma^{\lceil \lambda_j n\rceil-\lceil \lambda_{j-1} n\rceil}}(\Sigma)$ be the unique $\lceil \lambda_j n\rceil-\lceil \lambda_{j-1} n\rceil$-th level Bernoulli measure satisfying 
\begin{eqnarray*}
\nu_n([\omega_{\lceil \lambda_{j-1} n\rceil+1}\cdots\omega_{\lceil \lambda_j n\rceil}]):=\begin{cases} \frac{1}{\#\A^j_n}\text{ if }(\omega_{\lceil \lambda_{j-1} n\rceil+1},\cdots,\omega_{\lceil \lambda_j n\rceil})=\kappa^{\tau}\text{ for some }\tau\in \A^j_n\\
0\text{ otherwise.}
\end{cases}
\end{eqnarray*}
It follows that
\begin{eqnarray*}
\nu_n([\tau_{\lceil \lambda_{j-1} n\rceil+1}\cdots\tau_{\lceil \lambda_j n\rceil}]):=\begin{cases} \frac{1}{\#\A^j_n}\text{  for  }(\tau_{\lceil \lambda_{j-1} n\rceil+1},\cdots,\tau_{\lceil \lambda_j n\rceil})=\tau\in \A^j_n\\
0\text{  otherwise.}
\end{cases}
\end{eqnarray*}
Thus, $h\left(\nu_n\circ \chi_j^{-1},\sigma_j^{\lceil \lambda_j n\rceil-\lceil \lambda_{j-1} n\rceil}\right)=\log \#\A^j_n$ (see \cite{Walters} 4.26). Let $\mu_n:=A_{\lceil \lambda_j n\rceil-\lceil \lambda_{j-1} n\rceil}(\nu_n)$. By Lemma \ref{averagemeasure} (i) each $\mu_n$ is ergodic, and by Lemma \ref{averagemeasure} (ii) we have
\begin{eqnarray} \label{A n j ent bound}
\frac{1}{n}\log \#\B_n&=&\frac{{\lceil \lambda_j n\rceil-\lceil \lambda_{j-1} n\rceil}}{n}h(\mu^j_n\circ \chi_j^{-1},\sigma_j).
\end{eqnarray}

By the definition of $\Omega(A,N,q)$ for all $\omega\in \Omega(A,N,q)$ and $n> N$ there exists $\alpha_n$ such that for all $l\leq q$ we have 
\begin{eqnarray}
\bigg|\sum_{r=0}^{n-1} \varphi_l(\sigma^r(\omega))-n\int \varphi_l d\alpha_n\bigg|<\frac{n}{q}
\end{eqnarray}
and hence for all $n>\lambda_1^{-1}N$ we have 
\begin{eqnarray}
\bigg|\sum_{r=\lceil \lambda_{j-1} n\rceil}^{\lceil \lambda_{j} n\rceil-1} \varphi_l(\sigma^r(\omega))-(\lceil \lambda_j n\rceil-\lceil \lambda_{j-1} n\rceil)\int \varphi_l d\rho_n\bigg|<\frac{2\lceil \lambda_j n\rceil-\lceil \lambda_{j-1} n\rceil}{q}\nonumber.
\end{eqnarray}
Thus, by equation (\ref{omega of tau def}), for all $\tau \in A_n^j$, all $\omega \in [\omega^{\tau}_{\lceil \lambda_{j-1} n\rceil+1}\cdots \omega^{\tau}_{\lceil \lambda_{j} n\rceil}]$ and all $l\leq q$, 
\begin{eqnarray}
\bigg|\sum_{r=0}^{\lceil \lambda_{j} n\rceil-\lceil \lambda_{j-1} n\rceil-1} \varphi_l(\sigma^r(\omega))-(\lceil \lambda_j n\rceil-\lceil \lambda_{j-1} n\rceil)\int \varphi_l d\rho_n\bigg|\\
<\frac{2\lceil \lambda_j n\rceil-\lceil \lambda_{j-1} n\rceil}{q}+\sum_{r=0}^{\lceil \lambda_{j} n\rceil-\lceil \lambda_{j-1} n\rceil-1}\var_r(\varphi_l)\nonumber.
\end{eqnarray}
Since each $\varphi_l$ is continuous, it follows that there exists some $M>\lambda_1^{-1}N$ such that for all $n\geq M$, all $\tau \in A_n^j$, all $\omega \in [\omega^{\tau}_{\lceil \lambda_{j-1} n\rceil+1}\cdots \omega^{\tau}_{\lceil \lambda_{j} n\rceil}]$ and all $l\leq q$,
\begin{eqnarray}
\bigg|\frac{1}{\lceil \lambda_{j} n\rceil-\lceil \lambda_{j-1} n\rceil}\sum_{r=0}^{\lceil \lambda_{j} n\rceil-\lceil \lambda_{j-1} n\rceil-1} \varphi_l(\sigma^r(\omega))-\int \varphi_l d\rho_n\bigg|\\
<\frac{3\lambda_j - \lambda_{j-1}}{q( \lambda_j - \lambda_{j-1} )}<\frac{1}{\lceil Lq\rceil}.\nonumber
\end{eqnarray}
Now since $\nu_n$ is supported on sets of the form $[\omega^{\tau}_{\lceil \lambda_{j-1} n\rceil+1}\cdots \omega^{\tau}_{\lceil \lambda_{j} n\rceil}]$ with $\tau \in A_n^j$ it follows that for all $n\geq M$ and all $l\leq q$, 
\begin{eqnarray}
\bigg|\int A_{\lceil \lambda_{j} n\rceil-\lceil \lambda_{j-1} n\rceil}(\varphi_l) d\nu_n-\int \varphi_l d\rho_n\bigg|<\frac{1}{\lceil Lq\rceil}.\nonumber
\end{eqnarray}
Thus, by Lemma \ref{averagemeasure} (iii) we have
\begin{eqnarray}
\bigg|\int \varphi_l d\mu_n -\int \varphi_l d\rho_n\bigg|<\frac{1}{\lceil Lq\rceil}.\nonumber
\end{eqnarray}
for all $n\geq M$ and $l \leq \lceil Lq \rceil$, $\mu_n \in U(A,\lceil Lq\rceil)$ and hence 
\begin{equation}
h(\mu^j_n\circ \chi_j^{-1},\sigma_j)\leq H^j(A,\lceil Lq \rceil).
\end{equation}
By equation \ref{A n j ent bound} this proves the lemma.
\end{proof}

\begin{lemma} \label{PULemma} For each $q \in \N$,
\begin{equation*}
\Dim\Pi(\Omega(A))\leq  \frac{H^1(\sigma,A,\lceil Lq \rceil)}{\log a_1} +\sum_{j=2}^d \left(\frac{1}{\log a_j}-\frac{1}{\log a_{j-1}}\right) H^j(\sigma,A, \lceil Lq \rceil).
\end{equation*}
\end{lemma}
\begin{proof}
Combining Lemma \ref{First Box Estimate} with Lemma \ref{Entropy Bound} we have
\begin{equation*}
\BDim{\Pi(\Omega(A,N,q))}\leq  \frac{H^1(\sigma,A,\lceil Lq \rceil)}{\log a_1} +\sum_{j=2}^d \left(\frac{1}{\log a_j}-\frac{1}{\log a_{j-1}}\right) H^j(\sigma,A, \lceil Lq \rceil).
\end{equation*}
for each $N\in\N$. Moreover since
\begin{equation*}
\Pi(\Omega(A))\subseteq \bigcup_{N\in\N}\Pi(\Omega(A,N,q))
\end{equation*}
we may apply Proposition \ref{BoxPacking} to prove the lemma.
\end{proof}

\begin{lemma}\label{EntropyLimit} For each $j=1,\cdots,d$, $\lim_{l\rightarrow \infty}H^j(A,l)=H^j(\sigma,A).$
\end{lemma}
\begin{proof}
Fix $j\in \{1,\cdots,d\}$. Clearly $H^j(\sigma,A)\leq \liminf_{l\rightarrow \infty}H^j(A,l).$ Now for each $l\in \N$ choose $\mu_l \in U(A,l)$ with $h_{\mu_l\circ\chi_j^{-1}}(\sigma_j)>H^j(A,l)-\frac{1}{l}$. Since $\M_{\sigma}(\Sigma)$ is compact we may take a weak $*$ limit $\mu_{\infty} \in \M_{\sigma}(\Sigma)$. It follows from the fact that $A$ is closed and $\mu_l \in U(A,l)$ for each $l\in \N$, that $\mu_{\infty}\in A$. Moreover, since entropy is upper semi-continuous (see \cite[Theorem 8.2]{Walters}) 
\begin{eqnarray}
H^j(\sigma,A)&\geq& h_{\mu_{\infty}\circ\chi_j^{-1}}(\sigma_j) \\&\geq& \limsup_{l\rightarrow \infty}h_{\mu_l\circ\chi_j^{-1}}(\sigma_j)\nonumber \\
&\geq& \limsup_{l\rightarrow \infty}H^j(A,l)\nonumber.
\end{eqnarray}
\end{proof}

To complete the proof we let $q\rightarrow \infty$ in Lemma \ref{PULemma}. Applying Lemma \ref{EntropyLimit} we have
\begin{equation}
\Dim\Pi(\Omega(A))\leq  \frac{H^1(\sigma,A)}{\log a_1} +\sum_{j=2}^d \left(\frac{1}{\log a_j}-\frac{1}{\log a_{j-1}}\right) H^j(\sigma,A).
\end{equation}
This completes the proof of Theorem \ref{Packing Main symbolic} and hence Theorems \ref{Packing level sets} and \ref{Packing Main}.

\section{The Shape of the Spectrum}
We now deduce several features of the shape of the packing spectrum.
\begin{corollary}\label{concave}
Let $\varphi: \Sigma \rightarrow \R$ be a continuous real valued potential which is not cohomologous to a constant. Then, the packing spectrum $\alpha\mapsto \Dim{\Pi(J_{\varphi}(\alpha))}$ is concave and continuous on the interval $A(\varphi)=[\alpha_{\min},\alpha_{\max}]$.
\end{corollary}
\begin{proof} By Theorem \ref{Packing Main symbolic} it suffices to show that for each $j=\{1,\cdots,d\}$, $H^j(\sigma, \varphi, \alpha)$ is concave and continuous. So fix $j \in \{1,\cdots, d\}$. It follows from Lemma \ref{EntropyLimit} that $H^j(\sigma,\varphi,\alpha)$ is upper semi-continuous. Moreover, $H^j(\sigma,\varphi,\alpha)$ is concave. Indeed given $\alpha_-,\alpha_+\in A(\varphi)$ and $\delta>0$ we may choose $\mu_-,\mu_+\in \M_{\sigma}(\Sigma)$ such that $\int \varphi d\mu_-=\alpha_-$ and $\int \varphi d\mu_+=\alpha_+$ $h(\mu_-\circ \chi_j^{-1},\sigma_j)>H^j(\sigma,\varphi,\alpha_-)-\delta$, $h(\mu_+\circ \chi_j^{-1},\sigma_j)>H^j(\sigma,\varphi,\alpha_+)-\delta$. For each $t\in(0,1)$ we let $\mu_t:=(1-t)\mu_+t\mu_+$ so that $\int \varphi d\mu_t=(1-t)\alpha_-+t\alpha_+$. Moreover, since the entropy map is affine (see \cite[Theorem 8.1]{Walters})
\begin{eqnarray*}
H^j(\sigma,\varphi,(1-t)\alpha_-+t\alpha_+)&\geq& h(\mu_t\circ \chi_j^{-1},\sigma_j)\\
&=& (1-t)h(\mu_-\circ \chi_j^{-1}),\sigma_j)+th(\mu_+\circ\chi_j^{-1},\sigma_j)\\
&\geq& (1-t)H^j(\sigma_j,\varphi,\alpha_-)+tH^j(\sigma_j,\varphi,\alpha_+)-\delta
.\end{eqnarray*}
Letting $\delta\rightarrow 0$ we see that $\alpha \mapsto H^j(\sigma,\varphi,\alpha)$ is concave and hence lower semi-continuous. 
\end{proof}
The following corollary gives a sufficient condition on $\varphi$ for the packing spectrum to be analytic. 
\begin{corollary}
Suppose there is some H\"{o}lder continuous potential $\tilde{\varphi}:\Sigma_d\rightarrow \R$ such that $\varphi=\tilde{\varphi}\circ \chi_d$. Then $\alpha\mapsto \Dim{E_{\varphi}(\alpha)}$ is strictly concave and real analytic on the interval $A(\varphi)$.
\end{corollary}
\begin{proof} Note that for each $j=1,\cdots, d$, the projection $\chi_d\circ\chi_j^{-1}:\Sigma_j\rightarrow \Sigma_d$ is a well defined Lipchitz function. Hence the real valued potential $\varphi_j:\Sigma_j\rightarrow \Sigma_d$, given by $\varphi_j:=\tilde{\varphi}\circ \chi_d\circ\chi_j^{-1}$ is H\"{o}lder continuous. It follows straightforwardly from $\varphi=\tilde{\varphi} \circ \chi_d$ that for each $j=1,\cdots, d$,
\begin{equation}\label{conditional entropy variational}
H^j(\sigma,\varphi,\alpha)= \sup \left\lbrace h_{\mu}(\sigma_j): \mu \in \M_{\sigma_j}(\Sigma_j), \hspace{1mm}\int \varphi_j d\mu=\alpha\right\rbrace.
\end{equation}
One can deduce from standard results that the right hand side of (\ref{conditional entropy variational}) is strictly concave and analytic. Since $\varphi$ is not cohomologous to a constant and $\chi_j\circ \sigma = \sigma_j \circ \chi_j$ it is clear that no $\varphi_j$ is cohomologous to a constant. Now fix $j \in \{1,\cdots, d\}$. By \cite[Theorem 1.28]{Bowen} it follows that, for each $j$, the Gibbs measure corresponding to $\varphi_j$ is not the measure of maximal entropy on $\Sigma_j$. Now for each $\alpha \in A(\varphi)$ consider the set
\begin{equation}
J_j(\alpha):=\left\lbrace \omega \in \Sigma_j: \lim_{n\rightarrow \infty} \frac{1}{n}\sum_{r=0}^n \varphi_j(\sigma^r\omega)=\alpha\right\rbrace,
\end{equation}
where $\Sigma_j$ is given the usual symbolic metric (see \cite[Chapter 1]{Bowen}).
By \cite[Theorem 6]{Barreira Saussol} $\dim J_j(\alpha)$ is equal to a constant multiple of the quantity on the right hand side of (\ref{conditional entropy variational}). Since the Gibbs measure corresponding to $\varphi_j$ is not the measure of maximal entropy it follows from \cite[Theorem 1]{Pesin Weiss Birkhoff} that $\alpha \mapsto \dim J_j(\alpha)$ is strictly concave and real analytic on $(\alpha_{\min}, \alpha_{\max})$. Thus the spectrum is strictly convex and real analytic on $(\alpha_{\min}, \alpha_{\max})$. By Lemma \ref{concave} the spectrum is continuous on $[\alpha_{\min}, \alpha_{\max}]$ and hence these properties extend to the full interval $[\alpha_{\min},\alpha_{\max}]$.
\end{proof}

For each $j=1,\cdots, d$ we let $\mathbf{b}_j$ denote the measure of maximal entropy on $\Sigma_j$. We conclude this section with a necessary and sufficient condition for the packing spectrum to attain the full packing dimension of the repeller. The proof is immediate from Theorem \ref{Packing Main symbolic}.
\begin{corollary} There exists some $\alpha\in A(\varphi)$ satisfying $\Dim E_{\varphi}(\alpha)=\Lambda$ if and only if $\int \varphi d \mu_1=\int \varphi d \mu_2=\cdots =\int \varphi d \mu_d$ for some $\mu_1,\cdots, \mu_d \in \M_{\sigma}(\Sigma)$ such that $\mu_j\circ \chi_j^{-1}=\mathbf{b}_j$ for each $j=1, \cdots, d$.
\end{corollary}

\section{Examples}
In this section we consider two simple examples exhibiting interesting features of the packing spectrum. 

As noted in the introduction the packing and Hausdorff spectra need not coincide. This raises the question of whether there are any real-valued potentials $\varphi:\Sigma\rightarrow \R$ supported on Bedford-McMullen repellers for which $\dim(\Lambda)<\Dim(\Lambda)$ and yet the Hausdorff and packing spectra for $\varphi$ coincide. Our first example shows that this can indeed be the case. One consequence of this is that $
\Dim E_{\varphi}(\alpha)=\dim E_{\varphi}(\alpha) \leq \dim \Lambda <\Dim \Lambda$
for all $\alpha\in[\alpha_{\min},\alpha_{\max}]$. So the packing spectrum need not attain the full packing dimension of the repeller at any point. This is in contrast to the situation for Hausdorff dimension where there is always some $\alpha\in[\alpha_{\min},\alpha_{\max}]$ for which $\dim(E_{\varphi}(\alpha))=\dim(\Lambda)$, namely $\alpha=\int \varphi d\mu_{*}$ where $\mu_{*}$ is an invariant measure of full dimension (see  \cite{Barral Mensi} \cite{Bedford}, \cite{McMullen}). 
 
\begin{example}
Take $a_1=3$, $a_2=2$ and $\D=\left\lbrace (0,0),(1,1),(2,0)\right\rbrace$ and $\varphi:\Lambda\rightarrow \R$ defined by 
\begin{equation*}
\varphi(\omega)=\begin{cases} 1\text{  if   }\omega_1=(1,1)\\
0\text{  if }\omega_1\neq(1,1).\\
\end{cases}
\end{equation*}
Then, $\dim(\Lambda)<\Dim(\Lambda)$.
However, for all $\alpha\in [0,1]$,
\begin{eqnarray*}
\Dim E_{\varphi\circ \Pi}(\alpha)=\dim E_{\varphi\circ \Pi}(\alpha)=\frac{-\alpha\log\alpha -(1-\alpha)\log (1-\alpha)}{\log 2}+(1-\alpha)\frac{\log 2}{\log 3}
.\end{eqnarray*}
\end{example}

\begin{proof}
$\dim(\Lambda)<\Dim(\Lambda)$  follows from Theorem \ref{KP result}. By considering the $((1-\alpha)/2,\alpha,(1-\alpha)/2)$-Bernoulli measure, it follows from Proposition \ref{Barral Feng} that
\begin{eqnarray*}
\dim E_{\varphi}(\alpha)\geq \frac{-\alpha\log\alpha -(1-\alpha)\log (1-\alpha)}{\log 2}+(1-\alpha)\frac{\log 2}{\log 3}
.\end{eqnarray*}
It is easy to see that 
\[\sup\left\lbrace \sum_{i=1}^3-p_i\log p_i: p_i\in [0,1], \sum_{i=1}^3p_i=1, p_1=\alpha\right\rbrace  =
-\alpha\log\alpha -(1-\alpha)\log \left(\frac{1-\alpha}{2}\right)\]
\[\sup \left\lbrace \sum_{i=1}^2-p_i\log p_i: p_i\in [0,1], \sum_{i=1}^2p_i=1, p_1=\alpha\right\rbrace=
-\alpha\log\alpha -(1-\alpha)\log (1-\alpha).\] 
Moreover, it follows from the fact that $\varphi$ is locally constant (ie. $\varphi(\omega')=\varphi(\omega)$ for all $\omega,\omega'\in\Sigma$ with $\omega_1'=\omega_1$) together with the Kolmogorov-Sinai Theorem that the suprema
\begin{eqnarray*}
\sup \left\lbrace h_{\mu}(\sigma): \mu \in \B_{\sigma}(\Sigma), \int \varphi d\mu =\alpha \right\rbrace\\
\sup \left\lbrace h_{\mu \circ \chi_2^{-1}}(\sigma_v): \mu \in \B_{\sigma}(\Sigma), \int \varphi d\mu =\alpha \right\rbrace
\end{eqnarray*} 
are both attained by Bernoulli measures. Thus, applying Theorem \ref{Packing Main symbolic} we have
\begin{eqnarray*}
\Dim E_{\varphi\circ \Pi}(\alpha)\leq \frac{-\alpha\log\alpha -(1-\alpha)\log (1-\alpha)}{\log 2}+(1-\alpha)\frac{\log 2}{\log 3}
.\end{eqnarray*}
\end{proof}

For our next example we have identical $a_1$, $a_2$ and $\D$, along with a potential $\varphi$ which is prima facie very close to our previous one. Indeed for $\alpha\geq \frac{1}{2}$ the spectra for the two examples coincide (see Figure 2). However our next example has a point of non-analyticity at $\alpha=\frac{1}{2}$ and for $\alpha < \frac{1}{2}$ the two packing spectra are very different. In particular, the packing spectrum attains the packing dimension of $\Lambda$ and so rises above the Hausdorff spectrum.

\begin{example}
Take $a=3$, $b=2$ and $\D=\left\lbrace (0,0),(1,1),(2,0)\right\rbrace$ and $\varphi:\Lambda\rightarrow \R$ defined by 
\begin{equation*}
\varphi(\omega)=\begin{cases} 1\text{  if   }\omega_1=(2,0)\\
0\text{  if }\omega_1\neq(2,0)\\
\end{cases}
\end{equation*}
\begin{eqnarray*}
\Dim E_{\varphi\circ \Pi}(\alpha)=\begin{cases}\frac{-\alpha\log\alpha -(1-\alpha)\log (1-\alpha)-\alpha\log2}{\log 3}+1\text{ for }\alpha\leq \frac{1}{2}\\
\frac{-\alpha\log\alpha -(1-\alpha)\log (1-\alpha)}{\log 2}+(1-\alpha)\frac{\log 2}{\log 3}\text{ for }\alpha>\frac{1}{2}.
\end{cases}
\end{eqnarray*}
Moreover, $\alpha\mapsto \Dim E_{\varphi}(\alpha)$ is non-analytic and attains the full packing dimension $\Dim(\Lambda)$ at its maximum.
\end{example}

\begin{proof}
Note that 
\begin{eqnarray*}
\sup\left\lbrace \sum_{i=1}^3-p_i\log p_i: p_i\in [0,1], \sum_{i=1}^3p_i=1, p_1=\alpha\right\rbrace  =
-\alpha\log\alpha -(1-\alpha)\log \left(\frac{1-\alpha}{2}\right)
\end{eqnarray*}
\begin{eqnarray*}
\sup \left\lbrace \sum_{i=1}^2-p_i\log p_i: p_i\in [0,1], \sum_{i=1}^2p_i=1, p_1\geq \alpha\right\rbrace=\begin{cases}\log 2 \text{ for }\alpha\leq \frac{1}{2}\\
-\alpha\log\alpha -(1-\alpha)\log (1-\alpha) \text{ for }\alpha> \frac{1}{2}\end{cases}
.\end{eqnarray*}
Moreover, since $\varphi$ is locally constant the following suprema are both attained by Bernoulli measures,
\begin{eqnarray*}
\sup \left\lbrace h_{\mu}(\sigma): \mu \in \B_{\sigma}(\Sigma), \int \varphi d\mu =\alpha \right\rbrace\\
\sup \left\lbrace h_{\mu \circ \chi_2^{-1}}(\sigma_v): \mu \in \B_{\sigma}(\Sigma), \int \varphi d\mu =\alpha \right\rbrace
.\end{eqnarray*} 
 Thus, applying Theorem \ref{Packing Main symbolic} we have
\begin{eqnarray*}
\Dim E_{\varphi\circ \Pi}(\alpha)=\begin{cases}\frac{-\alpha\log\alpha -(1-\alpha)\log (1-\alpha)-\alpha\log2}{\log 3}+1\text{ for }\alpha\leq \frac{1}{2}\\
\frac{-\alpha\log\alpha -(1-\alpha)\log (1-\alpha)}{\log 2}+(1-\alpha)\frac{\log 2}{\log 3}\text{ for }\alpha>\frac{1}{2}
\end{cases}
.\end{eqnarray*}
Consequently, $\alpha\mapsto \Dim E_{\varphi}(\alpha)$ is non-analytic. This follows from the fact that the functions 
\begin{eqnarray*}
\alpha\mapsto \frac{-\alpha\log\alpha -(1-\alpha)\log (1-\alpha)-\alpha\log2}{\log 3}+1\\
\alpha\mapsto \frac{-\alpha\log\alpha -(1-\alpha)\log (1-\alpha)}{\log 2}+(1-\alpha)\frac{\log 2}{\log 3}
\end{eqnarray*}
 have distinct second derivatives at $\alpha=\frac{1}{2}$. It also follows from our expression for $\Dim E_{\varphi}(\alpha)$, together with \ref{KP result}, that the full packing dimension is attained at $\alpha=\frac{1}{3}$.
\end{proof}

\begin{figure}
\includegraphics[width=80mm]{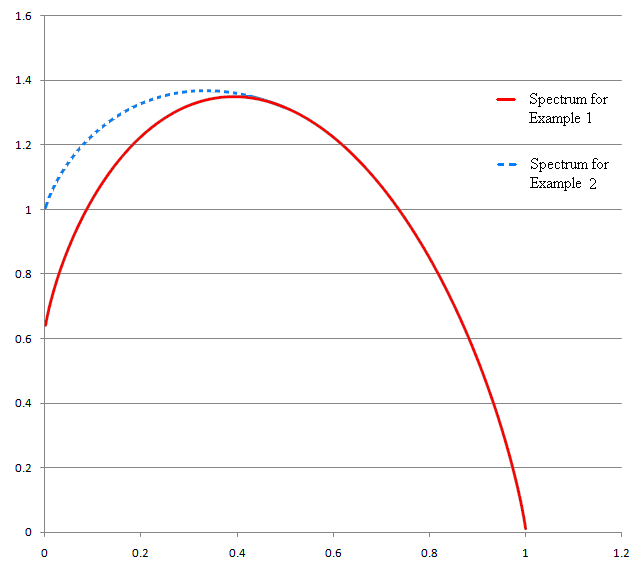}
\caption{The packing spectra for examples $1$ and $2$. The spectrum in Example $1$ is given by the red line. Example $2$ has a second order phase transition at $\alpha=\frac{1}{2}$. For $\alpha\leq \frac{1}{2}$ the spectrum in Example 2 is given by the broken blue line. For $\alpha > \frac{1}{2}$ the two spectra coincide and are both given by the red line.}
\end{figure}

\section{Generalisations and open questions}\label{Generalisations}

In this section we note some Corollaries to Theorem \ref{Packing Main} and \ref{Packing Main symbolic}. The first concerns sets of divergent points. Usually one considers sets of points for which the Birkhoff average converges to a given value. However, given any non-empty closed convex subset $A \subseteq A(\varphi)$ one may consider the set $E_{\varphi}(A)$ of points $x \in \Lambda$ for which the set of accumulation points for the sequence $(A_n(\varphi)(x))_{n\in\N}$ is equal to $A$. In the conformal setting both $\dim E_{\varphi}(A)$ and $\Dim E_{\varphi}(A)$ have been well studied in a series of papers due to Olsen and Winter \cite{Olsen Multifractal 1}, \cite{Olsen Multifractal 2}, \cite{Olsen Multifractal 3}, \cite{Olsen Multifractal 4}. This follows work by Barreira and Schmeling \cite{Nontypical1} showing that, given finitely many continuous potentials $\varphi_1,\cdots, \varphi_N$ on a conformal repeller $\Lambda$ for which each $A(\varphi_i)$ consists of at least two points, the Hausdorff dimension (and hence packing dimension) of the set of all points for which none of the Birkhoff averages for $\varphi_1,\cdots, \varphi_N$ converge is of full Hausdorff dimension. We note that \cite[Theorem 4.3]{Olsen Multifractal 4} implies that the set of points $x \in \Lambda$ for which the Birkhoff average $(A_n(\varphi)(x))_{n\in\N}$ does not converge for any continuous potential $\varphi$ for which $A(\varphi)$ consists of at least two points again, has full Hausdorff dimension. By a similar argument, along with some ideas from \cite{Kenyon Peres}, one can extend this result to self-affine Sierpi\'{n}ski sponges. 

One application of Theorem \ref{Packing Main symbolic} is to determine the packing dimension of the sets $E_{\varphi}(A)$ for self-affine Sierpi\'{n}ski sponges. We define, for $k=1,\cdots, d$,
\begin{eqnarray}
H^k(f,\varphi,A)&:=&\sup \left\lbrace h_{\mu\circ \pi_k^{-1}}(f_k): \mu \in \M(f,\Lambda), \int \varphi d\mu \in A \right\rbrace.
\end{eqnarray}

\begin{theorem}\label{Packing Divergence}
Let $\Lambda$ be a self-affine Sierpi\'{n}ski sponge. Let $\varphi:\Lambda \rightarrow \R^N$ be some continuous potential. Then given any non-empty closed convex subset $A \subseteq A(\varphi)$ we have
\[\Dim{E_{\varphi}(A)}=  \frac{H^1(f,\varphi,A)}{\log a_1} +\sum_{k=2}^d \left(\frac{1}{\log a_k}-\frac{1}{\log a_{k-1}}\right) H^k(f,\varphi,A).\]
\end{theorem}
In contrast very little is known concerning the Hausdorff dimension of $E_{\varphi}(A)$ for self-affine $\Lambda$, aside from the special case where $A$ is a singleton, and it would be very interesting to see if one could obtain a formula for $\dim E_{\varphi}(A)$ for arbitrary non-empty closed convex subsets of $\M(\Lambda,f)$.

Theorem \ref{Packing Main symbolic} also implies the some results concerning the packing spectrum for the local dimension of a Bernoulli measure on a self-affine Sierpi\'{n}ski sponge. To each Bernoulli measure $\mu$ on $\Sigma$ we associate the corresponding probability vector $(p_{i_1\cdots i_d})_{(i_1,\cdots, i_d)\in \D}$ in the usual way. Given $t\in \{1\,\cdots,d\}$ we let $p_{i_1\cdots i_t}$ denote the sum of all $p_{j_1\cdots j_d}$ for which $(j_1,\cdots, j_t)=(i_1,\cdots, i_t)$.

For $j=1, \cdots, d$ we define a potential $P_j:\Sigma \rightarrow \R$ by 
\begin{eqnarray}
P_j(\omega):= \begin{cases} \frac{\log p_{\chi_j(\omega_{1})}/ p_{\chi_{j+1}(\omega_{1})}}{\log a_j} \text{ if }j\neq d\\
\frac{\log p_{\chi_d(\omega_{1})}}{\log a_d}\text{ if }j=d.\end{cases}
\end{eqnarray}
Clearly $\var_1(P_j)=0$ for each $j$ and as such $P_j$ is continuous. Let $P:\Sigma \rightarrow \R^d$ denote the potential $\omega \mapsto (P_j(\omega))_{j=1}^d$. We shall assume the Very Strong Separation Condition (see \cite[Condition (II)]{Olsen Affine Pointwise}).

\begin{lemma}\label{Olsen symb point dim} Suppose that $\mu$ is a Bernoulli measure on a self-affine Sierpi\'{n}ski sponge which satisfies the Very Strong Separation Condition. Then for all $x=\Pi(\omega) \in \Lambda$ we have
\begin{equation}
\lim_{r\rightarrow \infty}\frac{\log \mu(B(x,r))}{\log r}=\sum_{j=1}^d \frac{1}{\lceil \lambda_j n \rceil}\sum_{l=0}^{\lceil \lambda_j n \rceil-1} P_j(\sigma^l\omega).
\end{equation}
\end{lemma}
\begin{proof} See \cite[Theorem 6.2.2]{Olsen Affine Pointwise}.
\end{proof}
Olsen \cite[Conjecture 4.1.7]{Olsen Affine Pointwise} conjectured that the packing spectrum of a Bernoulli measure on a self-affine Sierpi\'{n}ski sponge is given by the Legendre transform of an certain auxiliary function (see \cite[Section 3.1]{Olsen Affine Pointwise} for details). In particular, this conjecture would imply that the packing spectrum for local dimension always peaks at the full packing dimension of the attractor $\Lambda$ (see \cite[Theorem 3.3.2 (ix)]{Olsen Affine Pointwise} and note that $\gamma(0)= \Dim \Lambda$ by Theorem \ref{KP result}). Theorem \ref{Packing Main symbolic} provides us with the following counterexample.

\begin{example}
Take $a_1=4$, $a_2=3$ and $\D=\left\lbrace (0,0),(2,2),(3,0)\right\rbrace$ and let $\nu$ be the Bernoulli measure obtained by taking $p_{00}=p_{30}=1/4$ and $p_{22}=1/2$. Let $\alpha_{\min}:=\log 2/\log 3$ and $\alpha_{\max}:= \log 2/\log 3+1/2$ and for all $\alpha \in [\alpha_{\min},\alpha_{\max}]$ we define $\rho(\alpha):=2 \left(\alpha - \log 2/\log 3\right).$
Then, $\dim(\Lambda)<\Dim(\Lambda)$.
However, for all $\alpha\in [\alpha_{\min},\alpha_{\max}]$,
\begin{eqnarray*}
\Dim D_{\nu}(\alpha)=\dim D_{\nu}(\alpha)= \frac{-\rho(\alpha)\log \rho(\alpha) -(1-\rho(\alpha))\log (1-\rho(\alpha))}{\log 4}+\frac{1}{2}(1-\rho(\alpha))
.\end{eqnarray*}
\end{example}

\begin{proof}
Theorem \ref{KP result} implies $\dim(\Lambda)<\Dim(\Lambda)$. Applying Lemma \ref{Olsen symb point dim} we see that for all $\alpha \in [\alpha_{\min},\alpha_{\max}]$ we see that $\Pi(\omega) \in D_{\nu}(\alpha)$ if and only if
\begin{equation}
\V(\omega) \subseteq \left\lbrace \mu \in \M_{\sigma}(\Sigma): \mu([(0,0)])+\mu([3,0])=\rho(\alpha)\right\rbrace.
\end{equation}
Now proceed as in Example 1.
\end{proof}
Theorem \ref{Packing Main symbolic} also implies the following lower bound for the packing spectrum for local dimension. 
\begin{prop}\label{local dim packing lower bound}
Suppose that $\mu$ is a Bernoulli measure on a self-affine Sierpi\'{n}ski sponge which satisfies the Very Strong Separation Condition. Then,  
\begin{equation*}
\Dim D_{\mu}(\alpha) \geq \sup \left\lbrace \frac{H^1(\sigma,\underline{P},\underline{\alpha})}{\log a_1} +\sum_{k=2}^d \left(\frac{1}{\log a_k}-\frac{1}{\log a_{k-1}}\right) H^k(f,\underline{P},\underline{\alpha})\right\rbrace,
\end{equation*}
where the supremum is taken over all $\underline{\alpha}=(\alpha_j)_{j=1}^d\in \R^d$ for which $\sum_{j=1}^d\alpha_j=\alpha$.
\end{prop}
\begin{proof}
It follows from Lemma \ref{Olsen symb point dim} that $E_{\underline{P}}(\underline{\alpha})\subseteq D_{\mu}(\alpha)$ for each $\underline{\alpha}=(\alpha_j)_{j=1}^d\in \R^d$ with $\sum_{j=1}^d\alpha_j=\alpha$. Consequently, the result follows from Theorem \ref{Packing Main symbolic}.
\end{proof}
For a rather limited class of Bernoulli measures we obtain an equality.
\begin{defn}
We say that a Bernoulli measure $\mu$ on a self-affine Sierpi\'{n}ski sponge is one dimensional if there exists some $k\in \{1,\cdots,d\}$ for which the probability vector $(p_{i_1\cdots i_d})_{(i_1,\cdots, i_d)\in \D}$ associated to $\mu$ satisfies $p_{i_1\cdots i_{d+1-q}}/p_{i_1\cdots i_{d-q}}=\#\eta_{q+1}(\D)/\#\eta_q(\D)$ for all $(i_1,\cdots, i_{d+1-q})\in \eta_{q+1}(\D)$ and all $q \in \{1,\cdots,d-1\}\backslash \{k\}$ and each $p_i=1/\#\eta_d(\D)$ for $i \in \eta_d(\D)$, provided $k\neq d$.
\end{defn}
Now if $\mu$ is a one dimensional Bernoulli measure on a self-affine Sierpi\'{n}ski sponge then for each $j\neq k$ $P_j$ will be equal to an explicit constant $c_j$. Let $\tilde{P}$ denote the potential $P_k+\sum_{j\neq k}c_j$. 
\begin{theorem}
Suppose that $\mu$ is a one dimensional Bernoulli measure on a self-affine Sierpi\'{n}ski sponge which satisfies the Very Strong Separation Condition. Then 
\begin{equation*}
\Dim D_{\mu}(\alpha)= \frac{H^1(\sigma,\tilde{P},\alpha)}{\log a_1} +\sum_{k=2}^d \left(\frac{1}{\log a_k}-\frac{1}{\log a_{k-1}}\right) H^k(f,\tilde{P},\alpha).
\end{equation*}
\end{theorem}
\begin{proof} It follows from Lemma \ref{Olsen symb point dim} that $D_{\mu}(\alpha)=E_{\tilde{P}}(\alpha)$. Hence, the result follows from Theorem \ref{Packing Main symbolic}.
\end{proof}
We emphasise that the class of one dimensional Bernoulli measures is really very limited and the techniques of this paper are insufficient for determining $\Dim D_{\mu}(\alpha)$ for more general classes of Bernoulli measures. The reason for this extra level of difficulty is that one is essentially dealing with a sum of Birkhoff averages taken at multiple time scales (see Lemma \ref{Olsen symb point dim}). It seems unlikely that the lower bound given in Proposition \ref{local dim packing lower bound} is optimal. As such it remains an open question to determine the packing spectrum for local dimension on a self-affine Sierpi\'{n}ski sponge.


\begin{thebibliography}{25}
\bibitem[BOS]{Snigreva Packing}
 I. Baek, L. Olsen, N. Snigireva, \textit{Divergence points of self-similar measures and packing dimension}.  Adv. Math. 214 (2007), no. 1, 267-287.
\bibitem[BF]{Barral Feng}
J. Barral and D. Feng, \textit{Weighted thermodynamic formalism and applications}, (2009). arXiv:0909.4247v1.
\bibitem[BM1]{Barral Mensi}
J. Barral and M. Mensi, \textit{Multifractal analysis of Birkhoff averages on `self-affine' symbolic spaces}. Nonlinearity 21 (2008), no. 10, 2409-2425. 
\bibitem[BM2]{Barral Mensi Pointwise}
J. Barral and M. Mensi, \textit{Gibbs measures on self-affine Sierpinski carpets and their singularity spectrum}. Ergodic Theory Dynam. Systems 27 (2007), no. 5, 1419–1443. 
\bibitem[BSa]{Barreira Saussol}
L. Barreira and B. Saussol, \textit{Variational principles and mixed multifractal spectra}. Trans. Amer. Math. Soc. 353 (2001), no. 10, 3919-3944.
\bibitem[BSch]{Nontypical1}
L. Barreira and J. Schmeling, \textit{Sets of non-typical points have full topological entropy and full Hausdorff dimension}. Israel J. Math. 116 (2000), 29-70.
\bibitem[Be]{Bedford}
T. Bedford, \textit{PhD Thesis: Crinkly curves, Markov partitions and box dimension of self-similar sets}. Ph.D. thesis, University of Warwick (1984).
\bibitem[Bo]{Bowen}
R. Bowen, \textit{Equilibrium states and the ergodic theory of Anosov diffeomorphisms}. Second revised edition. Lecture Notes in Mathematics, 470. Springer-Verlag, Berlin, (2008). viii+75 pp. ISBN: 978-3-540-77605-5.
\bibitem[Ed1]{Edgar1}
G. A. Edgar,\textit{Measure, topology, and fractal geometry}. Second edition. Undergraduate Texts in Mathematics. Springer, New York, (2008). xvi+268 pp. ISBN: 978-0-387-74748-4.
\bibitem[Ed2]{Edgar2}
G. A. Edgar, \textit{Integral, probability, and fractal measures}. Springer-Verlag, New York, (1998). x+286 pp. ISBN: 0-387-98205-1
\bibitem[Fa1]{Falconer singular value}
K. Falconer, \textit{The Hausdorff dimension of self-affine fractals}. Math. Proc. Cambridge Philos. Soc. 103 (1988), no. 2, 339–350. 
\bibitem[Fa2]{Falconer}
K. Falconer, \textit{Fractal geometry. Mathematical foundations and applications}. Second edition. John Wiley $\&$ Sons, Inc., Hoboken, NJ, (2003). xxviii+337 pp. ISBN: 0-470-84861-8. 
\bibitem[Fa3]{Falconer Techniques}
K. Falconer, \textit{Techniques in fractal geometry}. John Wiley $\&$ Sons, Ltd., Chichester, 1997. xviii+256 pp. ISBN: 0-471-95724-0.
\bibitem[FFW]{Fan Feng Wu}
A. Fan, D. Feng and J. Wu, \textit{Recurrence, dimension and entropy}.  J. London Math. Soc. (2) 64 (2001), no. 1, 229--244. 
\bibitem[FLW]{Feng Lau Wu}
D. Feng K. Lau and J. Wu, \textit{Ergodic Limits on the Conformal Repellers}. Adv. Math. 169 (2002), no. 1, 58–91. 
\bibitem[GR]{Gelfert Rams}
K. Gelfert and M. Rams, \textit{The Lyapunov spectrum of some parabolic systems}. Ergodic Theory Dynam. Systems 29 (2009), no. 3, 919--940. 
\bibitem[JJOP]{non-uniformly hyperbolic}
A. Johansson, T. Jordan, A. Oberg, and M. Pollicott, \textit{Multifractal analysis of non-uniformly
hyperbolic systems}. Israel J. Math. 177 (2010), 125–144. 
\bibitem[JR]{Jordan Rams}
T. Jordan, M. Rams, \textit{ Multifractal analysis for Bedford-McMullen carpets}. Math. Proc. Camb. Phil. Soc. (2011), 147-156.
\bibitem[JS]{Jordan Simon}
T. Jordan, K. Simon, \textit{Multifractal Analysis of Birkhoff Averages for some Self-Affine IFS}. Dyn. Syst. 22 (2007), no. 4, 469–483.
\bibitem[Ki]{King}
J. King,\textit{ The Singularity Spectrum for General Sierpinski Carpets}. Adv. Math. 116 (1995), no. 1, 1-11.
\bibitem[KP]{Kenyon Peres}
R. Kenyon and Y. Peres, \textit{Measures of full dimension on affine-invariant sets}, Ergodic Theory Dynam. Systems 16 (1996), no. 2, 307–323.
\bibitem[Mat]{Mattila}
P. Mattila, \textit{Geometry of Sets and Measures in Euclidean Spaces: Fractals and rectifiability}. Cambridge Studies in Advanced Mathematics, 44. Cambridge University Press, Cambridge, (1995). xii+343 pp. ISBN: 0-521-46576-1; 0-521-65595-1.
\bibitem[McM]{McMullen}
C. McMullen, \textit{The Hausdorff Dimension of General Sierpinski Carpets}. Nagoya Math. J. 96 (1984), 1–9. 
\bibitem[Ni]{Nielsen}
O. Nielsen, \textit{The Hausdorff and packing dimensions of some sets related to Sierpi\'{n}ski carpets}. Canad. J. Math. 51 (1999), no. 5, 1073–1088. 
\bibitem[Ol1]{Olsen Multifractal 1}
L. Olsen, \textit{Multifractal analysis of divergence points of deformed measure theoretical Birkhoff averages}. J. Math. Pures Appl. (9) 82 (2003), no. 12, 1591–1649. 
\bibitem[Ol2]{Olsen Multifractal 3}
L. Olsen, \textit{Multifractal analysis of divergence points of deformed measure theoretical Birkhoff averages III.} Aequationes Math. 71 (2006), no. 1-2, 29–53.
\bibitem[Ol3]{Olsen Multifractal 4}
L. Olsen, \textit{Multifractal analysis of divergence points of deformed measure theoretical Birkhoff averages. IV: Divergence points and packing dimension}. Bull. Sci. Math. 132 (2008), no. 8, 650–678.
\bibitem[Ol4]{Olsen Affine Pointwise}
L. Olsen, \textit{Self-affine multifractal Sierpi\'{n}ski sponges in $\R^d$}. Pacific J. Math. 183 (1998), no. 1, 143–199.
\bibitem[OlWi]{Olsen Multifractal 2}
L. Olsen, S. Winter, \textit{Multifractal analysis of divergence points of deformed measure theoretical Birkhoff averages II: Non-linearity, divergence points and Banach space valued spectra}. Bull. Sci. Math. 131 (2007), no. 6, 518–558. 
\bibitem[Per]{Peres packing measure} 
Y. Peres, \textit{The packing measure of self-affine carpets}. Math. Proc. Camb. Phil. Soc. 115 (1994), no. 3, 437–450.
\bibitem[Pes]{Pesin}
Y. Pesin, \textit{Dimension Theory in Dynamical Systems. Contemporary views and applications}. Contemporary views and applications. Chicago Lectures in Mathematics. University of Chicago Press, Chicago, IL, 1997. xii+304 pp. ISBN: 0-226-66221-7; 0-226-66222-5.
\bibitem[PW1]{Pesin Weiss Birkhoff}
Y. Pesin, H. Weiss, \textit{The Multifractal Analysis of Birkhoff Averages and Large Deviations}. Global analysis of dynamical systems, 419–431, Inst. Phys., Bristol, 2001. 
\bibitem[PW2]{Pesin Weiss MA Equilibrium}
Y. Pesin, H. Weiss, \textit{The Multifractal Analysis of Gibbs Measures: Motivation, mathematical foundation and examples}.
Chaos 7 (1997), no. 1, 89–106. 
\bibitem[Wa]{Walters}
P. Walters, \textit{An Introduction to Ergodic Theory}. Graduate Texts in Mathematics, 79. Springer-Verlag, New York-Berlin, (1982). ix+250 pp. ISBN: 0-387-90599-5.
\end{thebibliography}
\end{document}